\documentclass[12pt]{article}
\usepackage{color}
\usepackage{amscd,amsfonts,amssymb,amscd,amsmath,latexsym}
\usepackage[hypertex]{hyperref}
\usepackage[all]{xy}
\usepackage{mathrsfs}

\title{Adams operations in smooth $K$-theory}
\author{Ulrich Bunke\thanks{NWF I - Mathematik,
Universit{\"a}t Regensburg,
93040 Regensburg,
GERMANY, ulrich.bunke@mathematik.uni-regensburg.de} }

\newtheorem{theorem}{Theorem}[section] 
\newtheorem{prop}[theorem]{Proposition}
\newtheorem{lem}[theorem]{Lemma}
\newtheorem{ddd}[theorem]{Definition}
\newtheorem{kor}[theorem]{Corollary}

\newcommand{\bR}{{\mathbf{R}}}

\newcommand{\bQ}{\mathbf{Q}}

\renewcommand{\P}{{\mathbb{P}}}

\newcommand{\Z}{\mathbb{Z}}

\newcommand{\bE}{{\bf E}}

\newcommand{\bU}{{\mathbf{U}}}

\newcommand{\proof}{{\it Proof.$\:\:\:\:$}}

\newcommand{\R}{\mathbb{R}}

\newcommand{\Q}{\mathbb{Q}}

\newcommand{\bK}{{\bf K}}

\newcommand{\bL}{{\bf L}}

\newcommand{\C}{\mathbb{C}}

\newcommand{\cE}{\mathcal{E}}

\newcommand{\cW}{\mathcal{W}}

\newcommand{\cK}{\mathcal{K}}

\newcommand{\Hom}{{\tt Hom}}

\newcommand{\im}{{\tt im}}

\newcommand{\inter}{{\tt int}}

\newcommand{\id}{{\tt id}}

\newcommand{\nat}{\mathbb{N}}

\newcommand{\cB}{\mathcal{B}}
\def\imath{{i}}

\def\hB{\hspace*{\fill}$\Box$ \newline\noindent}

\newcommand{\ind}{{\tt index}}

\def\hB{\hspace*{\fill}$\Box$ \\[0.5cm]\noindent}
\newcommand{\cL}{\mathcal{L}}
 
\newcommand{\cQ}{\mathcal{Q}}

\newcommand{\bH}{\mathbf{H}}
\newcommand{\bW}{\mathbf{W}}
\newcommand{\pr}{{\tt pr}}
\newcommand{\bX}{\mathbf{X}}
\newcommand{\bY}{\mathbf{Y}}

\newcommand{\ch}{{\mathbf{ch}}}

\newcommand{\bV}{\mathbf{V}}

\newcommand{\hA}{\hat{\mathbf{A}}}

\newcommand{\Bott}{{\mathtt{Bott}}}
\newcommand{\Thom}{{\mathtt{Thom}}}

\begin{document}\maketitle
\begin{abstract}
We show that the Adams operation $\Psi^k$, $k\in \{-1,0,1,2,\dots\}$, in complex $K$-theory lifts to an operation $\hat \Psi^k$ in smooth $K$-theory.
If $V\to X$ is a $K$-oriented vector bundle with Thom isomorphism $\Thom_V$, then
there is a characteristic class $\rho^k(V)\in K[\frac{1}{k}]^0(X)$ such that
$\Psi^k(\Thom_V(x))=\Thom_V(\rho^k(V)\cup \Psi^k(x))$ in $K[\frac{1}{k}](X)$ for all  $x\in K(X)$. 
We lift this class to a $\hat K^0(\dots)[\frac{1}{k}]$-valued characteristic class for real vector bundles with geometric $Spin^c$-structures.

If $\pi:E\to B$ is a $K$-oriented proper submersion, then for all $x\in K(X)$ we have 
$\Psi^k(\pi_!(x))=\pi_!(\rho^k(N)\cup \Psi^k(x))$ in $K[\frac{1}{k}](B)$, where $N\to E$ is the stable $K$-oriented normal bundle of $\pi$. To a smooth $K$-orientation $o_\pi$ of $\pi$  we associate a class $\hat \rho^k(o_\pi)\in \hat K^0(E)[\frac{1}{k}]$ refining $\rho^k(N)$. Our main theorem states that if $B$ is compact, then 
$\hat \Psi^k(\hat \pi_!(\hat x))=\hat \pi(\hat \rho^k(o_\pi)\cup \hat \Psi^k(\hat x))$ in $\hat K(B)[\frac{1}{k}]$ for all $\hat x\in \hat K(E)$.
We apply this result to the $e$-invariant of bundles of framed manifolds and $\rho$-invariants of flat vector bundles.

\end{abstract}

\tableofcontents

\section{Introduction}

The formalism of smooth extensions of  generalized cohomology theories is  designed to capture secondary invariants in topology, global analysis and geometry. The first example was the smooth extension of ordinary cohomology introduced by Cheeger-Simons \cite{MR827262}. Among other applications it was used to construct secondary characteristic classes for flat vector bundles.

Motivated by  applications in mathematical physics, in particular string theory, smooth extensions of other generalised cohomology theories, in particular of   $K$-theory, have been considered e.g in \cite{moore-2000-0005},  \cite{freed-2000-0005}  \cite{freed-2000}, \cite{szabo-2007}.  The existence of smooth extensions of  generalised cohomology theories has been shown in Hopkins-Singer \cite{MR2192936}. Axioms and uniqueness results have been discussed in \cite{MR2365651}, \cite{bunke-2009}.  In \cite{bunke-2009} we have shown that there is, up to unique isomorphism, a unique smooth extension of complex $K$-theory. 

An important tool for the construction of primary and secondary invariants is the integration or push-forward map for suitably oriented maps. The integration for smooth extensions of generalised  cohomology theories has been discussed in \cite{MR2192936}.  The notion of a smooth orientation of a submersion has been formalised for bordism theories in \cite{bunke-2007}, and in \cite{bunke-20071} for complex $K$-theory.

In general, Riemann-Roch type index theorems are assertions about the compatibility of natural operations between cohomology theories and the push-forward. In the prototypical example
its states the compatibility of the Chern character
$$\ch:K\to HP\Q$$
from $K$-theory to periodic rational cohomology with the push-forward along a $K$-oriented proper submersion $\pi:E\to B$ between smooth manifolds:
\begin{equation}\label{uiddwqdqwdqwd}
\xymatrix{K^*(E)\ar[d]^{\pi_!^K}\ar[r]^{\ch}&HP\Q^*(E)\ar[d]^{\pi_!^{HP\Q}(\hA^c(T^v\pi)\cup\dots)}\\
K^{*-n}(B)\ar[r]^\ch&HP\Q^{*-n}(B)}\ 
\end{equation}
Here $n:=\dim(E)-\dim(B)$ is the dimension of the fibres of $\pi$ and  $\hA^c(T^v\pi)\in HP\Q^0(E)$ is the $Spin^c$-generalisation  of the $\hA$-genus (see  \cite[Def. 3.3]{bunke-20071}) of the vertical bundle $T^v\pi:=\ker(d\pi)$ which has a $Spin^c$-structure by the $K$-orientation. The maps
$\pi^K_!$ and $\pi^{HP\Q}_!$ are the integration maps in the corresponding cohomology theories\footnote{In the main body of the paper we will omit the superscripts since we only consider
integration in $K$-theory}.

The prototypical result  for  the smooth extensions shown in  \cite[Thm. 6.19]{bunke-20071} states, that 
if $\pi$ is smoothly $K$-oriented, then the  smooth refinement of the diagram (\ref{uiddwqdqwdqwd})
\begin{equation}\label{uiddwqdqwdqwd1}
\xymatrix{\hat K^*(E)\ar[d]^{\hat \pi_!^K}\ar[r]^{\hat \ch}&\widehat{HP\Q}^*(E)\ar[d]^{\hat \pi_!^{HP\Q}(\hat\hA^c(o_\pi)\cup\dots)}\\
\hat K^{*-n}(B)\ar[r]^{\hat \ch}&\widehat{HP\Q}^{*-n}(B)}
\end{equation}
commutes, too.
Here $\hat K$ and $\widehat{HP\Q}$ denote the smooth extensions of complex $K$-theory and periodic rational cohomology theory, $\hat \ch:\hat K\to \widehat{HP\Q}$ is the smooth lift of the Chern character, and $\hat\hA^c(o_\pi)\in \widehat{HP\Q}^0(E)$
is the smooth refinement $\hA^c(o_\pi)$ determined by the smooth $K$-orientation $o_\pi$ of $\pi$.

In the present paper instead of the Chern character we consider the Adams operation
$$\Psi^k:K[\frac{1}{k}]^*(X)\to K[\frac{1}{k}]^*(X)$$  for $k\in \{-1\}\cup \nat$. In this case the  Riemann-Roch type theorem states that
\begin{equation}\label{uiddwqdqwdqwd2}
\xymatrix{K[\frac{1}{k}]^*(E)\ar[d]^{\pi_!^K}\ar[r]^{\Psi^k}&K[\frac{1}{k}]^*(E)\ar[d]^{\pi_!^K(\rho^k(T^v\pi)^{-1}\cup\dots)}\\
K[\frac{1}{k}]^{*-n}(B)\ar[r]^{\Psi^k}&K[\frac{1}{k}]^{*-n}(B)}\ 
\end{equation}
commutes, where $\rho^k(T^v\pi)\in K[\frac{1}{k}]^0(E)$ is an invertible $K$-theoretic characteristic class
of the $Spin^c$-bundle $T^v\pi$ (see Section \ref{uuiqwduqwdwqd}) below.

The main results of the present paper are the following three theorems:
\begin{theorem}[Theorem \ref{adams}]
There exists a  {\rm natural}  lift
 of the Adams operation to a natural transformation $\hat \Psi^k:\hat K(\dots)[\frac{1}{k}]\to \hat K(\dots)[\frac{1}{k}]$ of functors on the category of compact manifolds.
\end{theorem}

\begin{theorem}[Definition \ref{iqwdwqdqwdqwdwd} \& Theorem \ref{uidqwdqwd}]
If $\pi:E\to B$ is a submersion with compact $E$ which is smoothly $K$-oriented by $o_\pi$, then there exists a {\em natural} smooth refinement $\hat \rho^k(o_\pi)\in \hat K^0(E)[\frac{1}{k}]$ of the class
$\rho^k(T^v\pi)^{-1}$.
\end{theorem}
For details, in particular for the meaning of the word \textit{natural}, we refer to the main body of the present paper. The analog of (\ref{uiddwqdqwdqwd1}) is given by the third theorem.
\begin{theorem}[Theorem \ref{uifqfefewf5454577}]\label{uiddqwdqwdqwdqwdd}
If $\pi:E\to B$ is a smoothly $K$-oriented proper submersion over a compact base, then
 the smooth refinement of (\ref{uiddwqdqwdqwd2})  commutes:
\begin{equation}\label{uiddwqdqwdqwd3}
\xymatrix{\hat K^*(E)[\frac{1}{k}]\ar[d]^{\hat \pi_!^K}\ar[r]^{\hat \Psi^k}&\hat K^*(E)[\frac{1}{k}]\ar[d]^{\hat \pi_!^K(\hat \rho^k(o_\pi)\cup\dots)}\\
\hat K^{*-n}(B)[\frac{1}{k}]\ar[r]^{\hat \Psi^k}&\hat K^{*-n}(B)[\frac{1}{k}]}\ .
\end{equation}
\end{theorem}

The Dirac operator model of smooth $K$-theory  \cite{bunke-20071} provides the link between
between the push-forward in smooth $K$-theory and spectral geometric invariants of families of Dirac operators. So in principle, the diagram (\ref{uiddwqdqwdqwd3}) can be interpreted as a relation between these invariants for different families of Dirac operators. 
We discuss this aspect in greater detail in Section \ref{sec4}.

In \cite{bunke-20071},  we have constructed a version of Adams $e$-invariant $e(\pi)\in  K\R/\Z^{-n-1}(B)$ for families of $n$-dimensional framed manifolds $\pi:E\to B$ using only elements of the formalism of smooth $K$-theory. 
As an immediate consequence of (\ref{uiddwqdqwdqwd3}) we show
in Theorem \ref{udqiwdqwdqwd} that
$$k^L(\hat \Psi^k -1)e(\pi)=0$$ for sufficiently large $L$.  
In the case that $B$ is a point these relations for all $k\in \nat\cup\{-1\}$ together  imply the well-known (in fact optimal) upper bound of the  range of the $e$-invariant
\cite{MR0198468}.

\section{Adams operations}\label{uuiqwduqwdwqd}

Complex $K$-theory $K$ is a generalised cohomology theory. If we invert a number $k\in \{-1\}\cup \nat$, then we obtain the generalised cohomology theory $K[\frac{1}{k}]$.
For a finite  $CW$-complex $X$  we have 
\begin{equation}\label{zzzze}
K[\frac{1}{k}]^*(X)\cong K^*(X)[\frac{1}{k}]\ .
\end{equation}

By the Landweber formalism \cite{MR0423332} complex $K$-theory is   associated to the multiplicative formal group law $$(x,y)\mapsto x+y+bxy$$ over the ring $K^*:=\Z[b,b^{-1}]$ generated by the Bott element $b$ with $\deg(b)=2$.  The cohomology  theory $K[\frac{1}{k}] $ is then given by the same law considered over $K[\frac{1}{k}]^*=\Z[\frac{1}{k}][b,b^{-1}]$. The diagram
$$\xymatrix{Z[\frac{1}{k}][b,b^{-1}][[x]]\ar[rr]_{\Psi_{FGL}^k}^{x\mapsto kx,b\mapsto k^{-1}b}&&Z[\frac{1}{k}][b,b^{-1}][[x]]\\
\Z[\frac{1}{k}][b,b^{-1}]\ar[u]\ar[rr]^{b\mapsto k^{-1}b}_{\Psi^k}&&\Z[\frac{1}{k}][b,b^{-1}]\ar[u]}$$
gives a morphism $\Psi_{FGL}^k$ of formal group laws over 
the morphism of rings $\Psi^k$.  It induces 
the Adams operation $\Psi^k$ which is a multiplicative cohomology operation of the generalised cohomology theory $K[\frac{1}{k}]$.

It is the stable version of the classical Adams operation 
$$\Psi^k:K^0\to K^0$$
which is already defined before inverting $k$. 
 If $L\to X$ is a one-dimensional complex vector
bundle over a finite $CW$-complex $X$ and $[L]\in K^0(X)$
denotes the corresponding $K$-theory class, then we have
\begin{equation}\label{qqq12}
 \Psi^k([L])=[L^k]
\end{equation}
in $K^0(X)$.

The Bott periodicity isomorphism
$\Bott:K[\frac{1}{k}]^{*}(X)\stackrel{\sim}{\to} K[\frac{1}{k}]^{*+2}(X)$
is given by multiplication with the Bott element  $b\in K^2$ so that
 the following diagram commutes for all $n\in \Z$:
\begin{equation}\label{qqq14}\xymatrix{K[\frac{1}{k}]^*(X)\ar[r]^{k^{-n}\Psi^k}\ar[d]^{\Bott^n}&K[\frac{1}{k}]^*(X)\ar[d]^{\Bott^n}\\ K[\frac{1}{k}]^{*+2n}(X)\ar[r]^{\Psi^k}&K[\frac{1}{k}]^{*+2n}(X)} \ .\end{equation}
The Adams operations  satisfy \begin{equation}\label{qqq15}\Psi^k\circ \Psi^l=\Psi^{kl}\end{equation} (here we invert $k$ and $l$).

We define the multiplicative cohomology operation $\Psi^k_H:HP\Q^*\to HP\Q^*$ on the periodic rational cohomology $HP\Q^*(X):=H^*(X;K^*_\Q)$
such that
$\Psi^k_H(x)=
x$ for $x\in H^{2n}(X;\Q)$ or $x\in H^{2n-1}(X;\Q)$ and
$\Psi_H^k(b)=k^{-1}b$.
Periodic rational cohomology is the natural target of the  Chern character $\ch:K^*\to HP\Q^{*}$, and  we have
 \begin{equation}\label{qqq18}
\ch\circ \Psi^k=\Psi^k_H \circ \ch\ .  \end{equation}

A real $n$-dimensional vector bundle $V\to X$ with a $Spin^c$-structure is $K$-oriented. We have a Thom isomorphism
$$\Thom_V:K[\frac{1}{k}]^*(X)\to \tilde K[\frac{1}{k}]^{*+n}(\Thom(V))\ ,$$
where $\tilde K[\frac{1}{k}](\Thom(V))$ denotes the reduced $K[\frac{1}{k}]$-theory of the Thom space of $V$. There exists a unique invertible characteristic class
$$\rho^k(V)\in K[\frac{1}{k}]^0(X)$$ (called the cannibalistic class in  \cite{MR0198468}) such that
\begin{equation}\label{stern}
\Psi^k(\Thom(x))=\Thom(\rho^k(V)\cup \Psi^k(x))\ ,\quad \forall x\in K[\frac{1}{k}]^*(X)\ .
\end{equation}
A $K$-orientation of a proper submersion $\pi:E\to B$
is determined by a $Spin^c$-structure of the vertical bundle $T^v\pi=\ker(d\pi)$. If $\pi$ is $K$-oriented, then we have an integration map
$$\pi_!:K[\frac{1}{k}]^*(E)\to K[\frac{1}{k}]^{*-n}(B)\ ,$$
where $n=\dim(E)-\dim(B)$ is the dimension of the fibres. The compatibility of the Adams operations and the integration is expressed by the identity
\begin{equation}\label{eq20}
 \Psi^k(\pi_!(x))=\pi_!(\rho^k(T^v\pi)^{-1}\cup \Psi^k(x))\ ,\quad \forall x\in K[\frac{1}{k}]^*(E)\ .
\end{equation}
It is an immediate consequence of the usual construction   of $\pi_!$ and (\ref{stern}).

\section{The lift of the Adams operations} 

We consider the smooth extension $(\hat K,R,I,a,\int)$
of complex $K$-theory \cite{bunke-20071}, \cite{bunke-2009} on the category of compact manifolds.  We restrict to  compact   manifolds since we will frequently use the isomorphism  (\ref{zzzze}).  In order to avoid this restriction   one could alternatively start with a smooth extension of $K[\frac{1}{k}]$.

In the present paper it is useful to keep track of degrees properly. So for the domain of $a$ and the target of $R$ we will take the periodic  differential forms $\Omega P^*(M):=\Omega^*(M,K^*_\R)$ with $K_\R^*:=K^*\otimes \R$, and the corresponding cohomology is the periodic de Rham cohomology $HP_{dR}(M,K^*_\R)\cong HP^*(M)$. We define the natural transformation $\Psi^k_\Omega:\Omega P^*(M)\to \Omega P^*(M)$ of ring-valued functors such that
$\Psi^k_\Omega(\omega)= \omega$ for $\omega\in \Omega^{*}(M)$, and $\Psi_\Omega^k(b)=k^{-1}b$. It induces a corresponding transformation $\Psi^k_H$ on the periodic cohomology.

\begin{theorem}\label{adams}
There exist a unique natural transformation
$\hat \Psi^k:\hat K(\dots)[\frac{1}{k}]\to \hat K(\dots)[\frac{1}{k}]$ of set-valued functors
on the category of compact manifolds such that
\begin{equation}\label{qr1}I\circ \hat \Psi^k=\Psi^k\circ I\ ,\quad R\circ \hat \Psi^k=\Psi^k_\Omega\circ R\ ,\end{equation}
and
\begin{equation}\label{qr2}\xymatrix{\hat K^0(S^1\times M)[\frac{1}{k}]\ar[r]^{\hat \Psi^k}\ar[d]^\int& \hat K^0(S^1\times M)[\frac{1}{k}]\ar[d]^\int\\\hat K^{-1}(M)[\frac{1}{k}]\ar[r]^{\hat \Psi^k}&\hat K^{-1}(M)[\frac{1}{k}]}\end{equation}
commutes. 
The transformation $\hat \Psi^k$ preserves the ring structure and satisfies \begin{equation}\label{qr3}\hat \Psi^k\circ \hat \Psi^l=\hat \Psi^{kl}\ .\end{equation}
\end{theorem}
\proof
We first show that there is a unique  natural transformation of set-valued functors
$$\hat \Psi^k:\hat K^0\to \hat K^0$$
which satisfies (\ref{qr1}). We then show that this transformation 
preserves the ring structure and satisfies (\ref{qr3}).
Finally we extend
$\hat \Psi^k$ to all degrees using Bott periodicity
and (\ref{qr2}). 

The space $\bK_0:= \Z\times BU$ represents the homotopy type of the classifying space of the functor $K^0$.
We choose by  \cite[Prop 2.1]{bs2009} a sequence of compact manifolds $(\cK_i)_{i\ge 0}$ together with maps $$x_i:\cK_i\to \bK_0\ ,\quad \kappa_i:\cK_i\to \cK_{i+1}$$ such that
\begin{enumerate}
\item $\cK_i$ is homotopy equivalent to an $i$-dimensional $CW$-complex, 
\item $\kappa_i:\cK_i\to \cK_{i+1}$ is an embedding of a closed submanifold,
\item $x_i:\cK_i\to \bK_0$ is $i$-connected,
\item  $x_{i+1}\circ \kappa_i=x_i$\ .
\end{enumerate}
Let $u\in K^0(\bK_0)$ be the universal class represented by the identity map $\bK_0\to \bK_0$. By \cite[Prop. 2.6]{bs2009} we can further choose a sequence $\hat u_i\in \hat K^0(\cK_i)$ such that
$I(\hat u_i)=x_i^*u$ and 
$\kappa_i^*\hat u_{i+1}=\hat u_i$ for all $i\ge 0$.
By \cite[Lem. 3.8]{bs2009} for  $k\ge 2i+2$  we have $H^{2j+1}(\cK_k,\R)=0$ for all $j\le i$.

The requirements
$I(\hat \psi^k(u_i))=\Psi^k (I(\hat u_i))\in K^0(\cK_i)$ and $R(\hat \psi^k(\hat u_i))=\Psi^k_\Omega(R(\hat u_i))$ fix a class $\hat \psi^k(\hat u_i)\in \hat K^0(\cK_i)$ uniquely up to 
elements of the form $a(\alpha)$ for $\alpha\in F^{\ge i}HP^{-1}(\cK_i)$, where
\begin{equation}\label{zw1}
 F^{\ge l}HP^{-1}(X):=\bigoplus_{2j+1\ge l} b^{-j-1}H^{2j+1}(X;\R)\subseteq HP^{-1}(X)\ .
\end{equation}

Let $M$ be a compact  manifold. In the following we construct a map 
$\hat \Psi^k:\hat K^0(M)\to \hat K^0(M)$. Let $\hat y\in \hat K^0(M)$ be given. 
Then we choose $i>\dim(M)$ and $f:M\to \cK_i$
such that $I(\hat y)=f^*x_i^*(u)$.
We further choose $\rho\in \Omega P^{-1}(M)$  such that

\begin{equation}\label{qr18}\hat y=f^*\hat u_i+a(\rho)\ .\end{equation}
We define
\begin{equation}\label{qr16}
\hat \Psi^k(\hat y):=f^*\psi^k(\hat u_i)+a(\Psi^k_\Omega(\rho ))\ . \end{equation}
By a direct calculation we verify that (\ref{qr1}) holds true:
$$I(\hat \Psi^k(\hat y))=\Psi^k(I(y))\ ,\quad R(\hat \Psi^k(\hat y))=\Psi^k_\Omega (R(\hat y))\ .$$
 \begin{lem} The right-hand side of (\ref{qr16})
does not depend on the choices. \end{lem} 
\proof If $\rho^\prime$ is a second choice for $\rho$ in (\ref{qr18}), then
$\rho^\prime-\rho=\ch(x)$ in $\Omega P^{-1}(M)/\im(d)$ for some $x\in K^{-1}(M)$. But then
$\Psi^k_\Omega(\rho)-\Psi^k_\Omega(\rho)=\Psi^k_H(\ch(x))=\ch(\Psi^k(x))$. This implies $a(\Psi^k_\Omega\rho)=a(\Psi^k_\Omega\rho^\prime)$.

Since we take $i$ sufficiently large it follows that
$f^*\psi^k(\hat u_i)$ is independent of the choice of the actual element $\psi^k(\hat u_i)$. 

We further can increase $i$ to $i+1$ by replacing $f$ by $\kappa_i\circ f$ without changing the right-hand side of (\ref{qr16}).

 Finally, any two choices of $f$ become homotopic after increasing $i$ sufficiently many times. Therefore let $f^\prime$ be another choice for $f$ in (\ref{qr18}) and $H:[0,1]\times M\to \cK_i$ be a homotopy from $f$ to $f^\prime$.
We choose $\tilde \rho\in \Omega P^{-1}([0,1]\times M)$ such that $\pr_M^*\hat y=H^*\hat u_i+a(\tilde \rho)$. In (\ref{qr18})
we can then take $\rho=\tilde \rho_{|\{0\}\times M}$ and $\rho^\prime:=\tilde \rho_{|\{1\}\times M}$, and by the homotopy formula (see \cite[(1)]{bs2009}) we have
$$a\left(\rho^\prime-\rho+\int_{[0,1]\times M/M} H^* R(\hat u_k)\right)=0\ ,$$
i.e. $$\rho^\prime-\rho+\int_{[0,1]\times M/M} H^* R(\hat u_k)=\ch(x)$$ for some $x\in K^{-1}(M)$.

We evaluate the difference of the right-hand sides of (\ref{qr16}) for the two choices and get
\begin{eqnarray*}
  \lefteqn{f^{\prime,*}\psi^k(\hat u_i)+a(\Psi^k_\Omega(\rho^\prime))-f^{ *}\psi^k(\hat u_i)-a(\Psi^k_\Omega(\rho)  )}&&\\&=& a\left(\int_{[0,1]\times M/M}R(H^*\psi^k(\hat u_k))+\Psi^k_\Omega(\rho^\prime)-\Psi^k_\Omega(\rho) \right) \\&=&
 a\left(\Psi^k_\Omega\left( \int_{[0,1]\times M/M} H^* R(\hat u_k)+ \rho^\prime- \rho\right)\right)  \\&=&a( \ch(\Psi^k(x)))\\&=&0\ ,
\end{eqnarray*}
where we use the homotopy formula in the first equality. \hB
We now have constructed for each manifold $M$ a map of
sets $\hat \Psi^k:\hat K^0(M)\to \hat K^0(M)$ satisfying (\ref{qr1}).
\begin{lem}
$\hat \Psi^k$ is a natural transformation of set-valued functors on the category of compact manifolds. It is the unique natural  transformation satisfying (\ref{qr1}). 
\end{lem}
\proof
Let $g:M^\prime\to M$ be a smooth map of compact manifolds. If we take in addition $i> \dim(M^\prime) $, then we can start the construction with 
$g^*\hat y=f^\prime \hat u_i+a(g^*\rho)$, the analog of (\ref{qr18}), where 
$f^\prime:=f\circ g$. With choice we have.  
$$\hat \Psi^k (g^*\hat y)= (f\circ g)^*\psi^k(\hat u_i)+a(\Psi^k_\Omega g^*\rho)
=g^*(f^*\psi^k(\hat u_i)+a(\Psi^k_\Omega(\rho)))=g^*\hat \Psi^k(\hat y)\ .$$

We now show uniqueness. If $\Psi^{\prime,k}:\hat K^0\to \hat K^0$ is another  natural transformation satisfying (\ref{qr1}), then
for $i>\dim(M)$ and $f:M\to \cK_i$ we have
\begin{equation}\label{qz1}
 \hat \Psi^{\prime,k} (f^*\hat u_i)=f^*\hat \Psi^{\prime,k}(\hat u_i)=f^*\psi^k(\hat u_i)=\hat \Psi^k(f^*\hat u_i)\ .
\end{equation}

For $\rho\in \Omega P^{-1}(M)$ we consider the class $\hat y:=f^*\hat u_i+a(\rho)$ and
$\hat x:=\pr_M^*f^*\hat u_i+a(t \pr_M^*\rho)\in \hat K^0([0,1]\times M)$, where $t$ is the coordinate of $[0,1]$.
The homotopy formula gives
\begin{eqnarray*}
\hat \Psi^{k,\prime}(\hat y)-\hat \Psi^{k,\prime}(f^*\hat u_i)
&=&a\left(\int_{[0,1]\times M/M]} R(\hat \Psi^{\prime,k}(\hat x))\right)\\
&=&a\left(\int_{[0,1]\times M/M]}\Psi^k_{\Omega} (R(\hat x))\right)\\
&=& \hat \Psi^{k }(\hat y)-\hat \Psi^{k }(f^*\hat u_i)\ .
\end{eqnarray*}
In view of (\ref{qz1}) we get
$\hat \Psi^{k,\prime}(\hat y)=\hat \Psi^{k }(\hat y)$.
\hB

\begin{lem} \label{lll3}
 $\hat \Psi^k:\hat K^0\to \hat K^0$ is a natural transformation of ring-valued functors and satisfies (\ref{qr3}).
\end{lem}
\proof
We first consider the additive structure.
Let $$\hat B(\hat x,\hat y):=\hat \Psi^k(\hat x+\hat y)-\hat \Psi^k(\hat x)-\hat \Psi^k(\hat y)\ .$$
Since $\hat \Psi^k$ is compatible with $I$ and $R$ we immediately see that
$\hat B$ takes values the subfunctor $H P^{-1}/\im(\ch)\subset \hat K^0$.
Furthermore, since by the explicit formula (\ref{qr16}) we have 
$\hat \Psi^k(\hat y+a(\rho))=\hat \Psi^k(\hat y)+\Psi^k(a(\rho))$, it follows  that
$\hat B$ factorises over a natural transformation
$$B:K^0\times K^0\to H P^{-1}/\im(\ch)\ .$$
The same argument as for Theorem \cite[3.6]{bs2009} shows that such a transformation vanishes.
This shows that $\hat \Psi^k$ preserves the additive structure.

In order to show that $\hat \Psi^k$ is multiplicative we argue similarly.
We consider
$$\hat E(\hat x,\hat y):=\hat \Psi^k(\hat x\cup \hat y)-\Psi^k(\hat x)\cup \Psi^k(\hat y)\ .$$ 
We again see that $\hat E$ factors over a transformation
$$E:K^0\times K^0\to HP^{-1}/\im(\ch)$$
which necessarily vanishes.

For the relation $\hat \Psi^l\circ \hat \Psi^k=\hat \Psi^{lk}$ we argue similarly using
$$\hat C(\hat x):=  \Psi^l\circ \hat \Psi^k(\hat x)-\hat \Psi^{lk}(\hat x)\ .$$
We again see that
$\hat C$ factors over a natural transformation
$$C:K^0\ \to H P^{-1}/\im(\ch)$$
which necessarily vanishes. \hB

\begin{lem}
There exists a unique natural transformation
$\hat \Psi^k:\hat K^{-1}\to \hat K^{-1}$ such that (\ref{qr2}) commutes. It further satisfies (\ref{qr1}) and (\ref{qr3}), and
for $\hat z\in \hat K^0(M)$ and $\hat y\in \hat K^{-1}(M)$ we have 
\begin{equation}\label{qz6}\hat \Psi^k(\hat z\cup \hat y)=\hat \Psi^k(\hat z)\cup\hat  \Psi^k(\hat y)\ .\end{equation}
\end{lem}
\proof
Let $\hat e\in \hat K^1(S^1)$ be characterised
by the properties $\int \hat e=1$ and that $R(\hat e)$ is rotation invariant.
For $\hat x\in \hat K^{-1}(M)$ by  (\ref{qr2}) we are forced to define
$$\hat \Psi^k(\hat x):=\int \hat \Psi^k(\hat e\times \hat x)\ .$$
This gives a natural transformation such that (\ref{qr2})  commutes. The relations
 (\ref{qr1}), (\ref{qr3}) and (\ref{qz6}) follow by direct calculations.
 \hB

By the relations (\ref{qr1}) and (\ref{qqq14}) we are forced to extend the transformation $\hat \Psi^k$ to all degrees by Bott periodicity, i.e. such that
\begin{equation}\label{qqq141}\xymatrix{\hat K^*(M)[\frac{1}{k}]\ar[r]^{k^{-n}\hat\Psi^k}\ar[d]^{\Bott^n}&\hat K^*(M)[\frac{1}{k}]\ar[d]^{\Bott^n}\\ \hat K^{*+2n}(M)[\frac{1}{k}]\ar[r]^{\hat \Psi^k}&\hat K^{*+2n}(M)[\frac{1}{k}]}\end{equation}
commutes. The relation (\ref{qr3}) holds true automatically.

In order to finish the proof we must show that
$\hat \Psi^k$ is multiplicative.
Let $\hat x,\hat y\in \hat K^*(M)$.
Then we can write
$\hat x=b^l \hat x_1$, $\hat y=b^n \hat y_1$, where
$\hat x_1,\hat y_1$ have degrees in $\{0,-1\}$. In this way by (\ref{qqq141}) we reduce the problem to 
the multiplicativity of $\hat \Psi^k$ in degree zero and (\ref{qz6}). 
This finishes the proof of Theorem \ref{adams}. \hB

An alternative way to construct the lift of the Adams operations would be to use the model \cite[Thm. 2.5]{bunke-2007}
$$\hat K(M):=\widehat{MU}(M)\otimes_{MU^*}K^*\ .$$ 
 together with \cite[Cor. 2.8]{bunke-2007}.

The following example shows that the lift of the Adams operations to smooth $K$-theory act on the classes of geometric line bundles in the expected way lifting (\ref{qqq12}). 
Let $\bL:=(L,h^L,\nabla^L)$ be a hermitian line bundle with connection over $M$.
It gives rise to a geometric family $\cL$ (see \cite[2.1.4]{bunke-2007}) and a smooth $K$-theory class $[\bL]:=[\cL,0]\in \hat K^0(M)$ in the model of smooth $K$-theory (compare \cite[Lemma 2.16]{bunke-2007}).
\begin{prop}\label{udiduqwdwqddwqdwqdw}
We have $\hat \Psi^k([\bL])=[\bL^k]$.
\end{prop}
\proof
We first consider the canonical bundle
$U\to \C\P^n$. Equipped with geometry  we get the geometric bundle
$\bU=(U,h^U,\nabla^U)$ and the class $[\bU]\in \hat K^0(\C\P^n)$.
Note that by a direct calculation
$$I(\hat \Psi^k([\bU]))=  I([\bU^k])\ ,\quad R(\hat \Psi^k[\bU])  = R([\bU^k])\ .$$
Since $HP^{-1}(\C\P^n)=0$ the class
$\hat \Psi^k([\bU])$ is uniquely determined by its curvature and its underlying topological $K$-theory class. This implies
\begin{equation}\label{qz2}\hat \Psi^k([\bU])=[\bU^k]\ .\end{equation}

In the general case there exists $n\ge 0$ and  $f:M\to \C\P^n$ such that
$L\cong f^*U$. We consider the bundle
$K:=\pr_M^*L\to [ 0,1]\times M$ with a geometry $\bK$, which coincides with
$f^*\bU$ on $\{0\}\times M$ and with $\bL$ on $\{1\}\times M$.
From the homotopy formula \cite[(1)]{bs2009} we get
\begin{eqnarray*}
\hat \Psi^k([\bL])-f^*\hat \Psi^k([\bU])&=&
a\left(\int_{[0,1]\times M/M} R(\hat \Psi^k([\bK]))\right)\\
&=&a\left(\int_{[0,1]\times M/M}   R([\bK^{k}])\right)\\
&=&[\bL^k]-f^*[\bU^k]\ .
\end{eqnarray*}
In view of (\ref{qz2}) this implies the assertion. \hB

\newcommand{\hHQ}{{\widehat{HP\Q}}}
Let $(\hHQ,R,I,a,\int)$ denote the smooth extension of
the periodic rational cohomology theory.
In \cite{bunke-20071}
we have constructed a lift of the Chern character to a natural transformation
$$\hat \ch:\hat K\to \hHQ$$
of ring-valued functors. We let
$\hat \Psi^k_H:\hHQ^*\to \hHQ^*$ denote the obvious lift of
$\Psi^k_H$ which multiplies  $b^n$ by  $k^{-n}$.
\begin{prop}\label{u9qwdqwdqwdqwd}
We have
$$\hat \ch\circ \hat \Psi^k=\hat \Psi^k_H\circ \hat \ch\ .$$
\end{prop}
\proof
We first consider the even case.
The difference
$$\hat D:=
\hat \ch\circ \hat \Psi^k-\hat \Psi^k_H  \circ \hat \ch$$
factors over a matural transformation
$$D:K^{0}\to HP^{-1}/HP\Q^{-1} \ .$$
Using that $\bK_0$ is an even space similar arguments as in the proof of Lemma \ref{lll3} 
show that $D=0$. The odd case follows from the compatibility of the Adams operations and the Chern character with integration.
\hB 

Let $M$ be a compact connected  manifold with base point $*\in M$. We consider the  multiplicative subgroups 
$$ U:=\{x\in   K^0(M)[\frac{1}{k}]\:|\: x_{|*}=1\}\ ,\quad \hat U:=\{\hat x\in \hat K^0(M)[\frac{1}{k}]\:|\:\hat x_{|*}=1\}$$
of $K^0(M)[\frac{1}{k}]$ and $\hat K^0(M)[\frac{1}{k}]$.
\begin{lem}\label{suskwswqoidd}
The maps
$$\Psi^k:U\to U\ ,\quad \hat \Psi^k:\hat U\to \hat U$$ are surjective. 
\end{lem}
\proof
We first consider the topological case.
 Let $F^{2n} K^0(M)[\frac{1}{k}]\subseteq K^0(M)[\frac{1}{k}]$ denote the $2n$'th step of the Atiyah-Hirzebruch filtration which is finite.
 The Atiyah-Hirzebruch filtration is compatible with the multiplication in the sense that
$$F^{2n} K(M)[\frac{1}{k}]\cup F^{2m} K(M)[\frac{1}{k}]\subseteq F^{2n+2m}K(M)[\frac{1}{k}]\ .$$
Let $x\in U$. Then we find inductively
approximations $z\in U$ such that $\Psi^k(z)-x\in F^{2n}K^0[\frac{1}{k}](M)$ as follows.
The first approximation is
$z=1$. Then
$\Psi^k(1)-x\in F^2 K^0(M)[\frac{1}{k}]$.
Assume that
$\Psi^k(z)-x=:d\in F^{2n-2}K^0(M)[\frac{1}{k}]$.
Then we take
$$z^\prime:=z-\frac{1}{k^n} d\ .$$
Then $\Psi^k(\frac{1}{k^n}d)-d\in F^{2n} K^0[\frac{1}{k}](M)$ 
and therefore
$\Psi^k(z^\prime)-x=\Psi^k(z)-\frac{1}{k^n}\Psi^k(d)-x\in  F^{2n} K^0(M)[\frac{1}{k}]$.

We now consider the smooth case.
Let $\hat x\in \hat U$. Then $I(\hat x)\in U$. We thus can choose
$z\in U$ such that $\Psi^k(z)=I(\hat x)$.
Let $\hat z$ be a smooth lift.
Then
$\hat \Psi^k(\hat z)-\hat x=a(\omega)$
for some $\omega\in \Omega P^{-1}(M)$.
We define
$\hat z^\prime:=\hat z-a((\Psi^k_\Omega)^{-1}(\omega))$.
Then
$\hat \Psi^k(\hat z^\prime)=\hat x$. \hB

\section{The characteristic class $\hat \rho^k$}

A $Spin^c$-structure $(P,\phi)$ on a real $n$-dimensional vector bundle $V\to M$ is a pair of a $Spin^c$-principal bundle
$P\to M$ together with an isomorphism
$\phi:P\times_{Spin^c(n)} \R^n\stackrel{\sim}{\to} V$, where $Spin^c(n)$ acts on $\R^n$ via the natural projection $Spin^c(n)\to SO(n)$. 
A $Spin^c$-vector bundle is $K$-oriented.
Let $\Thom_V:K[\frac{1}{k}]^*(M)\to \tilde K[\frac{1}{k}]^{*+n}(\Thom(V))$ denote the
Thom isomorphism. Then there is a characteristic class
$\rho^k(V)\in K[\frac{1}{k}]^0(M)$ of $Spin^c$-vector bundles uniquely characterised by  the relation
$$\Psi^k(\Thom_V(x))=\Thom_V(\rho^k(V)\cup \Psi^k(x))\ ,\quad  \forall x\in K[\frac{1}{k}]^*(M)\ ,$$
see \cite{MR0198468}.
We consider the characteristic class
$\hA^c(V)\in HP^0(M)$ of $Spin^c$-bundles. A definition is given in
\cite[Def.3.3]{bunke-20071}. In  the present paper we modify this definition by inserting suitable powers of the Bott element $b\in K^*$ in order to shift all homogenous components to degree zero. The $K$-theoretic characteristic class $\rho^k$
has the following properties:
\begin{enumerate}
\item $\rho^k(V\oplus W)=\rho^k(V)\cup \rho^k(W)$,
\item $\rho^k(M\times \R^n)=1$,
\item $\rho^k(V)$ is a unit in $K[\frac{1}{k}]^0(M)$,
\item $\ch(\rho^k(V))=\frac{\Psi^k_H( \hA^c(V))}{\hA^c(V)}$.
\end{enumerate}
A geometric $Spin^c$-structure on a vector bundle $V\to M$ is a triple
$\bV=(P,\phi,\tilde \nabla)$, where $(P,\phi)$ is a $Spin^c$-structure on $V$
  and $\tilde\nabla$ is a
connection on $P$. By Chern-Weil theory we can define
a closed differential form
$\hA^c(\bV)\in \Omega P^0(M)$ which represents the cohomology class $\hA^c(V)$, see \cite[Def. 3.3]{bunke-20071}

There is a natural definition of the sum $\bV\oplus \bW$ of two geometric $Spin^c$-vector bundles.
\newcommand{\vect}{\mathrm{Vect^{Spin^c}}}
We consider the  contravariant functor
$$\vect:\mathtt{smooth \:\:compact\:\:manifolds}\longrightarrow \mathtt{semi\:\: groups}$$
which associates to a compact smooth manifold $M$  the semigroup $\vect(M)$ of geometric $Spin^c$ vector bundles.
\begin{theorem}\label{uidqwdqwd}
There exists a unique natural transformation of set-valued functors
$$\hat \rho^k:\vect\to \hat K^0(\dots )[\frac{1}{k}]$$
such that
\begin{eqnarray}
I\circ \hat \rho^k&=&\rho^k \label{udidqwdwqdwd1}\\
 R(\hat \rho^k(\bW))&=&\frac{\Psi^k_\Omega(\hA^c(\bW)) }{\hA^c(\bW)}\ ,\quad \bW\in \vect(M)\label{udidqwdwqdwd}\ .\end{eqnarray}
 
It in addition satisfies
 \begin{eqnarray}
 \hat \rho^k(\bW\oplus \bW^\prime)&=&\hat \rho^k(\bW)\cup \hat \rho^k(\bW^\prime)\ ,\quad   \bW,\bW^\prime\in \vect(M)\label{ui11}\\
 \hat \rho^k(\bW) &&\mbox{is a unit in $\hat K^0(M)[\frac{1}{k}]$ for $\bW\in \vect(M)$}.
\end{eqnarray}  
\end{theorem}
\proof

Let us fix $n\in \nat_0$. We first construct $\hat \rho^k$ on the subfunctor
$\vect_n$ of $n$-dimensional geometric $Spin^c$-vector bundles.
Note that the classifying space $BSpin^c(n)$ of $Spin^c(n)$ is a simply-connected rationally even space. Using
   \cite[Prop 2.1]{bs2009} we choose a sequence of compact manifolds $(\cB_i)_{i\ge 0}$ together with maps $$x_i:\cB_i\to BSpin^c(n)\ ,\quad \kappa_i:\cB_i\to \cB_{i+1}$$ such that
\begin{enumerate}
\item $\cB_i$ is homotopy equivalent to an $i$-dimensional $CW$-complex, 
\item $\kappa_i:\cB_i\to \cB_{i+1}$ is an embedding of a closed submanifold,
\item $x_i:\cB_i\to BSpin^c(n)$ is $i$-connected,
\item  $x_{i+1}\circ \kappa_i=x_i$\ .
\end{enumerate}
Let $\xi_n\to BSpin^c(n)$ be the universal $n$-dimensional $Spin^c$ vector bundle. Let $E_i:=x_i^*\xi_n$. 
Note that $\kappa_i^* E_{i+1}\cong E_i$.  Since $\kappa_i:\cB_i\to \cB_{i+1}$ is an embedding of a closed submanifold we can 
inductively choose geometric refinements   $\bE_i\in \vect_n(\cB_i)$  such that 
$\kappa_i^*\bE_{i+1}\cong \bE_i$.  For each $i$ we
 choose a class $\hat r^k(\bE_i)\in \hat K^0(\cB_i)[\frac{1}{k}]$ such that
$$I(\hat r^k(\bE_i))=\rho^k(E_i)\ ,\quad  R(\hat r^k(\bE_i))=\frac{\Psi^k_\Omega(\hA^c(\bE_i))}{\hA^c(\bE_i)}\ .$$ This element is uniquely determined up to elements in
$a(F^{\ge i}HP^{-1}(\cB_i))$, see (\ref{zw1}). 
 
Let now $M$ be a compact   manifold and $\bW\in \vect_n(M)$.
Then we choose $i\in \nat $ such that $i>\dim(M)$ and there exists $f:M\to \cB_i$  and an isomorphism $f^*E_i\cong W$. We are forced to define
$\hat \rho^k(f^*\bE_i)=f^*\hat r^k(\bE_i)$. This class
is independent of the choice of $\hat r^k(\bE_i)$.

We now consider the bundle $V:=\pr_M^*W\to [0,1]\times M$. We choose a geometry $\bV\in \vect_n([0,1]\times M)$ which coincides on $\{0\}\times M$ with $f^*\bE_i$, and on $\{1\}\times M$ with $\bW$. 
By the homotopy formula we are forced to define
\begin{equation}\label{t12}\hat \rho^k(\bW):=\hat \rho^k(f^*\bE_i)+a\left(\int_{[0,1]\times M/M} R(\hat \rho^k(\bV))\right)=f^*\hat r^k(\bE_i)+a\left(\int_{[0,1]\times M/M} \frac{\Psi^k_\Omega (\hA^c(\bV))}{\hA^c(\bV)}\right).\end{equation}
By a straightforward calculation we verify that (\ref{udidqwdwqdwd1}) and (\ref{udidqwdwqdwd}) hold true. 
\begin{lem}
 $\hat \rho^k:\vect_n\to \hat K^0(\dots )[\frac{1}{k}]$ is a well-defined natural transformation.
\end{lem}
\proof
If we choose a second geometry $\bV^\prime\in \vect_n([0,1]\times M)$ interpolating between
$\bW$ and $f^*\bE_i$, then
the difference
$$\int_{[0,1]\times M/M} \frac{\Psi^k_\Omega (\hA^c(\bV))}{\hA^c(\bV)} -\int_{[0,1]\times M/M} \frac{\Psi^k_\Omega (\hA^c(\bV^\prime))}{\hA^c(\bV^\prime)}$$
is exact. This follows by Stokes' theorem from the fact, that 
the geometries $\bV$ and $\bV^\prime$ can again be connected by geometric
bundle over $[0,1]^2\times M$.
We conclude that $\hat \rho^k(\bW)$ does not depend on the choice of 
$\bV$.
 
If we increase $i$ by one and set $f^\prime:=\kappa_i\circ f$, then
we get $f^{\prime,*}\bE_{i+1}\cong f^*\bE_i$, $f^*\hat r^k(\bE_i)=f^{\prime,*} \hat r^k(\bE_{i+1})$,
 and therefore the same result for
$\hat \rho^k(\bW)$. 

Any two choices of maps $f:M\to \cB_i$ become homotopic after increasing $i$ sufficiently many times.
 If $f$ and $f^\prime$ are homotopic by a homotopy 
$H:[0,1]\times M\to \cB_i$,
then we apply the construction for $\pr_M^*\bW$ and $H$.
In this way we get a class
$\hat\rho^k(\pr_M^*\bW)$. Note that $$R(\hat \rho^k(\pr_M^*\bW))=\pr_M^*\frac{\Psi^k_\Omega (\hA^c(\bW))}{\hA^c(\bW)}$$ has no $dt$-component, where $t$ is the coordinate of $[0,1]$.
It thus follows from the homotopy formula that
$$\hat\rho^{k,\prime}(\bW)=\hat \rho^k(\pr_M^*\bW)_{|\{1\}\times M}=\hat \rho^k(\pr_M^*\bW)_{|\{0\}\times M}
 =\hat \rho^k(\bW)\ .$$

 Finally we  verify that $\hat \rho^k:\vect_n\to \hat K^0(\dots)[\frac{1}{k}]$ is a natural transformation. Let $g:M^\prime\to M$ be a smooth map. Then in the definition of $\hat \rho^k(g^*\bW)$ we can take
$f^\prime:=f\circ g$ and $\bV^\prime:=(\id_{[0,1]}\times g)^*\bV$. With these choices we have
$$g^*\hat \rho^k(\bW)=\hat \rho^k(g^*\bW)\ .$$
\hB

\begin{lem}
The relation (\ref{ui11}) holds true. 
\end{lem}
\proof
We consider the transformation
$$B:\vect_n\times \vect_m\to \hat K^0(\dots)[\frac{1}{k}]$$ given by
$$B(\bW,\bW^\prime):=\hat \rho^k(\bW\oplus \bW^\prime)-\hat \rho^k(\bW)\cup \hat \rho^k(\bW^\prime)\ .$$
We must show that $B=0$.
By construction we have
$$R\circ B=0\ ,\quad I\circ B=0\ .$$
Therefore
$B$ takes values in the homotopy invariant subfunctor
$$H P^{-1}(\dots)/\im(\ch)[\frac{1}{k}]\subset  \hat K^0(\dots)[\frac{1}{k}]\ .$$
Since it is homotopy invariant it is clear that
$B(\bW,\bW^\prime)\in H P^{-1}(M)/\im(\ch)[\frac{1}{k}]$ only depends on the underlying topological $Spin^c$-bundles $W$ and $W^\prime$ over $M$.
 There exists $j>\dim(M)$ such that
  can find  maps $f,f^\prime:M\to \cB_j$ such that
$W\cong f^*E_j$ and $W^\prime\cong f^{\prime,*} E_j$.
We set $F:=f\times f^\prime$, $V:=\pr_1^*E_j$, and $V^\prime:=\pr_2^*E_j$, where  $\pr_k:\cB_j\times \cB_j\to \cB_j$, $k=1,2$ are the projections. Then we get 
$$B(W,W^\prime)=B(F^*V,F^*V^\prime)=F^*B(V,V^\prime)\ .$$
Now $B(V,V^\prime)\in F^{\ge j}  H P^{-1}(\cB_j\times \cB_j)/\im(\ch)[\frac{1}{k}]$, and therefore
 $F^*B(V,V^\prime)=0$.
This shows that $B(W,W^\prime)=0$.
\hB

\begin{lem} For $\bW\in \vect(M)$ 
 the class $\hat \rho^k(\bW)\in \hat K^0(M)[\frac{1}{k}]$ is a unit. 
\end{lem}
\proof
We write
$\hat \rho^k(\bW)=1+(\hat \rho^k(\bW)-1)$. It suffices to show that
$\hat \rho^k(\bW)-1$ is nilpotent.
First of all, since $$R(\hat \rho^k(\bW))=1+\mbox{\em higher \:order\: forms}$$ we see that
$R(\hat \rho^k(\bW)-1)$ is nilpotent.
Similarly, the restriction of $$I(\hat \rho^k(\bW)-1)=\rho^k(W)-1$$ to a point vanishes.
Therefore $\rho^k(W)-1$ belongs to a lower step of the Atiyah-Hirzebruch filtration
of $K^0(M)[\frac{1}{k}]$ and is therefore nilpotent.

We conclude that for some large $l\in \nat $
$$(\hat \rho^k(\bW)-1)^l\in  H P^{-1}(M)/\im(\ch)[\frac{1}{k}]\ .$$
But then
$(\hat \rho^k(\bW)-1)^{2l}=0$.
\hB 
This finishes the proof of Theorem \ref{uidqwdqwd}. \hB

In view of the homomorphism $Spin(n)\to Spin^c(n)$ a $Spin$-structure on $V$  naturally induces a $Spin^c$-structure. We define the notion of a  geometric $Spin$-bundle in a similar manner as the notion of a geometric $Spin^c$-bundle. Notice that a geometric $Spin$-bundle $\check{\bW}$ gives rise to a geometric
$Spin^c$-bundle $\bW$.  We have $\hA^c(\bW)=\hA(\check \bW)$,  and this form is invariant under $\Psi^{-1}_\Omega$.

\begin{lem}
 If $\check{\bW}$ is a geometric $Spin$-bundle, then
$\hat \rho^{-1}(\bW)=1$.
\end{lem}
\proof
If $\check{\bW}$ is an $n$-dimensional geometric $Spin$-bundle on a compact manifold, then there exist another geometric $Spin$-bundle $\check{\bE}$ on a manifold
$B$ which has no real cohomology in odd degree below $\dim(M)+1$, and a map $f:M\to B$ such that $f^*E\cong W$. In fact, for $B$ we can take some approximation of a finite skeleton of $BSpin(n)$.
Notice that
$\hat \rho^{-1}(\check{\bE})-1\in F^{\ge \dim(M)+1} H P^{-1}(B)/\im(\ch)$.
Therefore $\hat \rho^{-1}(f^*\bE)=1$.
We now consider a geometric $Spin$-bundle $\check{\bV}$
over $[0,1]\times M$ which connects $\check{\bW}$ and
$f^*\bE$. Since $\hA(\check{\bV})=\hA^c(\bV)$ is invariant under $\Psi^{-1}_\Omega$  we have
$R(\hat \rho^{-1}(\bV))=1$. If we use these observations in (\ref{t12}) we get $\hat \rho^{-1}(\bW)=1$.
\hB

Let us now consider a  submersion $\pi:E\to B$ from a compact manifold $E$. We assume that the vertical bundle $T^v\pi:=\ker(d\pi)$ has a $Spin^c$-structure. It induces a $K$-orientation of $\pi$.
Recall that a smooth $K$-orientation $o$ (\cite[Def. 3.5]{bunke-20071}) of $f$ is represented  by a tuple
$(g^{T^v\pi},T^h\pi,\tilde \nabla,\sigma)$, where $g^{T^v\pi}$ is a vertical metric,
$T^h\pi$ is a horizontal distribution, $\tilde \nabla$ is a $Spin^c$-extension of the Levi-Civita connection $\nabla^{T^v\pi}$ on
$T^v\pi$ (induced by $g^{T^v\pi}$ and $T^h\pi$), and $\sigma \in \Omega P^{-1}(E)/\im(d)$ (here we again use the modified definition based on the insertion powers of the Bott element in order to shift all forms to degree $-1$).
The smooth $K$-orientation in particular induces a geometric $Spin^c$-structure $\mathbf{T^v\pi}$ on $T^v\pi$.
The curvature of the $K$-orientation is defined by 
$$R(o):=\hA^c(\tilde \nabla)-d\sigma\ ,$$
where we write  $\hA^c(\tilde \nabla)$
instead of $\hA^c(\mathbf{T^v\pi})$.

 
We consider the bundle $\bar \pi=\id\times \pi:[0,1]\times E\to [0,1]\times B$ and choose a representative of a  smooth $K$-orientation
$\bar o$ which interpolates from
$(g^{T^v\pi},T^h\pi,\tilde \nabla,0)$ to
$(g^{T^v\pi},T^h\pi,\tilde \nabla,\sigma)$.
\begin{ddd}\label{iqwdwqdqwdqwdwd} We define
$$\hat \rho^k(o):=\hat \rho^k(\mathbf{T^v\pi})^{-1}+a\left(\int_{[0,1]\times E/E}  \frac{\Psi^k_\Omega(R(\bar o))}{R(\bar o)}\right)\in \hat K^0(E)[\frac{1}{k}]\ .$$
\end{ddd}
Note that two tuples $(g^{T^v\pi},T^h\pi,\tilde \nabla,\sigma)$ and $(g^{T^v\pi,\prime},T^{h,\prime}\pi,\tilde \nabla^\prime,\sigma^\prime)$
represent the same smooth $K$-orientation if the underlying $Spin^c$-structures are isomorphic and
$\sigma^\prime-\sigma=\tilde \hA^c(\tilde\nabla^\prime,\tilde\nabla)$.
Here $\tilde \hA^c(\tilde\nabla^\prime,\tilde\nabla)$ is the transgression form defined in \cite[Def. 3.4]{bunke-20071} (again shifted to degree $-1$).

\begin{prop} \label{udqidqwdqwdqwdqwdd}
The class $\hat \rho^k(o)$ is   independent  of the choice of the representative of $o$.
\end{prop}
\proof
We choose a representative of a smooth $K$-orientation $\bar{\bar o}$ of $[0,1]\times E\to [0,1]\times B$ which interpolates from $(g^{T^v\pi,\prime},T^{h,\prime}\pi,\tilde \nabla^\prime,\sigma^\prime)$ to $(g^{T^v\pi},T^{h}\pi,\tilde \nabla,\sigma)$ .
Furthermore we let $\bar o^\prime$ be the smooth $K$-orientation
 of $[0,1]\times E\to [0,1]\times B$ obtained
by concatenating $\bar {\bar o}$ with $\bar o$.
With these choices we get
$$
\hat \rho^{k,\prime}(o)-\hat \rho^k(o)= a\left(
 \int_{[0,1]\times E/E} \frac{\Psi^k_\Omega(R(\bar{\bar o}))}{R(\bar{\bar o})}\right) \ .
 $$
In order to go further we adopt a very special choice for $\bar{\bar \sigma}$:
$$\bar{\bar\sigma}(t):=\pr_E^*\sigma^\prime+\int_{[0,t]\times E/E} \hA^c(\bar{\bar \nabla})\ .$$
Indeed, 
$$\bar{\bar\sigma}(1)=\pr_E^*\sigma^\prime+\hA^c(\tilde \nabla,\tilde \nabla^\prime)=\pr_E^*\sigma\ .$$
Then we have
$d\bar{\bar\sigma}=dt\wedge i_{\partial_t} \hA^c(\bar{\bar \nabla})+\theta$
with $i_{\partial_t}\theta=0$.
It follows that
$i_{\partial_t}R(\bar{\bar o})=0$ and therefore
$$ \int_{[0,1]\times E/E} \frac{\Psi^k_\Omega(R(\bar{\bar o}))}{R(\bar{\bar o})}  =0\ .$$
 \hB

We now consider an iterated bundle
$$\xymatrix{W\ar@/^1cm/@{..>}[rr]^r\ar[r]^p&E\ar[r]^q&B}$$
of compact manifolds.
We assume that
$T^vp$ and $T^vq$ are equipped with $Spin^c$-structures.
We choose smooth $K$-orientations
$o_p$ and $o_q$ lifting these $Spin^c$-structures. Then there is an induced $Spin^c$-structure on
$T^vr\cong T^vp\oplus p^*T^vq$ and a smooth orientation
$o_r=o_q\circ o_p$ (see \cite[Def. 3.21]{bunke-20071})

\begin{prop}
We have
$\hat \rho^k(o_p\circ o_q)=\hat \rho^k(o_p)\cup p^*\hat \rho^k(q)$.
\end{prop}
\proof
We consider the difference
$$\Delta(o_p,o_q):=\hat \rho^k(o_q\circ o_p)-\hat \rho^k(o_p)\cup p^*\hat \rho^k(q)\ .$$
We first check by a direct calculation that
$$ I(\Delta(o_p,o_q))= 0\ , \quad 
 R(\Delta(o_p,o_q))=0\ .$$  
It follows that
$\Delta(o_p,o_q)\in H P^{-1}(W)/\im(\ch)$.
Since two choices of smooth $K$-orientations refining a fixed underlying topological $K$-orientation can be connected by a path
it  follows by homotopy invariance that
$\Delta(o_p,o_q)$ only depends on the topological $Spin^c$-structures of $T^vp$ and $T^vq$.
Let us now recall the construction of
$o_r$.  We take $o_p:=(g^{T^vp},T^hp,\tilde\nabla^{T^vp},0)$ and
$o_q:=(g^{T^vq},T^hq,\tilde\nabla^{T^vq},0)$.
Then for $\lambda>0$ we get an induced metric
$g^{T^vr}_\lambda=\lambda^2g^{T^vp}\oplus p^*g^{T^vq}$ and splitting
$T^hr$. As explained in \cite[3.3.1]{bunke-20071} we also get an induced connection
$\tilde \nabla^{T^vr}_\lambda$ which has a limit $\tilde \nabla^{adia}:=\lim_{\lambda\to 0}\tilde \nabla^{T^vr}_\lambda$.
The composition of the smooth $K$-orientations is represented by (see \cite[Def. 3.21]{bunke-2007})
$$o_q\circ o_p=(g^{T^vr}_\lambda,T^hr,\tilde \nabla^{T^vr}_\lambda,-\tilde \hA^c(\tilde \nabla^{adia},\tilde \nabla^{T^vr}_\lambda))\ .$$
In \cite[(20)]{bunke-2007} we have observed that
$\tilde\nabla^{adia}=\tilde \nabla^{T^vp}\oplus \tilde p^*\nabla^{T^vq}$ under
the decomposition $\mathbf{T^vr}\cong \mathbf{T^vp}\oplus p^*\mathbf{T^vq}$ of vector bundles with geometric $Spin^c$-structures. 
We have by Definition \ref{iqwdwqdqwdqwdwd} and Proposition \ref{udqidqwdqwdqwdqwdd} that
\begin{eqnarray*}\hat \rho^k(o_q\circ o_p)&=&\hat \rho^k(\mathbf{T^vr})^{-1}+a(T(\lambda))\\
&=&\hat \rho^k(\mathbf{T^vp}\oplus p^*\mathbf{T^vq})^{-1}+a(T^\prime (\lambda))\\
&=&\hat \rho^k(\mathbf{T^vp})^{-1}\cup p^*\hat \rho^k(\mathbf{ T^vq})^{-1}+a(T^\prime(\lambda))\\
&=&\hat \rho^k(o_p)\cup p^*\hat \rho^k(o_q)+a(T^\prime(\lambda))\ ,
\end{eqnarray*}
where the forms $T(\lambda)$ and $T^\prime(\lambda)$ depend on the difference
of $\tilde \nabla^{adia}$ and $\tilde \nabla^{T^vr}_\lambda$ and vanish as $\lambda\to 0$.
 We now take the limit $\lambda\to 0$ and get
$$\hat \rho^k(o_q\circ o_p)=\hat \rho^k(o_p)\cup p^*\hat \rho^k(o_q)\ .$$
\hB

We now consider a cartesian diagram
$$\xymatrix{F\ar[r]^g\ar[d]^q&E\ar[d]^p\\ A\ar[r]^f&B}\ .$$
A smooth $K$-orientation $o_p$ of $p$ induces a smooth $K$-orientation $o_q$ of $q$.
\begin{lem}\label{uiqwdwqdwqdwqdwqd}
In this situation we have
$$\hat \rho^k(o_q)=g^*\hat \rho^k(o_p)\ .$$
\end{lem}
 \proof
The geometric $Spin^c$-structure on $T^vq$ induced by $o_q$ is
$\mathbf{T^vq}\cong g^*\mathbf{T^vp}$.  If we choose $\bar o_q:=(\id_{[0,1]}\times g)^*\bar o_p$
in the construction of $\hat \rho^k(o_q)$, then we immediately get 
$\hat \rho^k(o_q)=g^*\hat \rho^k(o_p)$ from  Definition \ref{iqwdwqdqwdqwdwd}. 
\hB

Consider a submersion $\pi:E\to B$ from a compact manifold $E$.
  In \cite[5.11]{bunke-2007} we have observed that 
a stable framing of  $T^v\pi$ provides a canonical $K$-orientation. 
Let us assume that $T^v\pi\oplus \underline{\R^N}_E$ is framed, where
$\underline{\R^N}_E:=E\times \R^N$ denotes the trivial $N$-dimensional real vector bundle over $E$. The associated smooth $K$-orientation is 
$$o_\pi:=(g^{T^v\pi},T^h\pi,\tilde \nabla^{T^v\pi}, \tilde \hA^c(\tilde \nabla^{T^v\pi}\oplus\nabla^{\underline{\R^N}_E } ,\tilde \nabla^{T^v\pi\oplus \underline{\R^N}_E, frame}))\ ,$$ where
$\tilde \nabla^{T^v\pi\oplus \underline{\R^N}_E,frame}$ is the connection  induced by the faming.

\begin{prop}\label{udiqwdwqdw33333}
If $o_\pi$ is induced by a stable framing of $T^v\pi$, then 
$\hat \rho^k(o_\pi)=1$.
\end{prop}
\proof
 We have by the homotopy formula
\begin{equation}\label{ui1}\hat \rho^k(\mathbf{T^v\pi})^{-1}= 1+a\left(\int_{[0,1]\times E/E} \frac{\Psi^k_\Omega (\hA^c(\bar \nabla))}{\hA^c(\bar \nabla)}\right) \ ,\end{equation}
 where $\bar \nabla$ is a family of $Spin^c$-connection interpolating from $\tilde \nabla^{T^v\pi\oplus \underline{\R^N}_E,frame}$  to $\tilde \nabla^{T^v\pi}\oplus \tilde \nabla^{\underline{\R^N}_E }$.
 Let $\bar \sigma$ be a form on $[0,1]\times E$ which interpolates from $\tilde \hA^c(\tilde \nabla^{T^v\pi}\oplus\tilde \nabla^{\underline{\R^N}_E } ,\tilde \nabla^{T^v\pi\oplus \underline{\R^N}_E, frame})$ to zero. This family induces a smooth
$K$-orientation $\bar o$ on $[0,1]\times E\to [0,1]\times B$ which interpolates
from $o_\pi$ to
$(g^{T^v\pi},T^h\pi,\tilde \nabla^{T^v\pi},0)$.
 Then we have  
\begin{equation}\label{ui2}\hat \rho^k(o_\pi)= \hat \rho^k(\mathbf{T^v\pi})^{-1}-
a\left(\int_{[0,1]\times E/E}\frac{\Psi^k_\Omega(R(\bar o))}{R(\bar o)}\right) \ .\end{equation}
Furthermore, we have 
$$\tilde \hA^c(\tilde \nabla^{T^v\pi}\oplus\tilde \nabla^{\underline{\R^N}_E } ,\tilde \nabla^{T^v\pi\oplus \underline{\R^N}_E, frame}) =\int_{[0,1]\times E/E} \hA^c(\bar \nabla)\ .$$
As in the proof of Proposition \ref{udqidqwdqwdqwdqwdd} we can choose
$$\bar \sigma(t):=\int_{[t,1]\times E/E}\hA^c(\bar \nabla)$$
so that
\begin{eqnarray*}d\bar \sigma(t)&=&-dt\wedge i_{\partial_t} \hA^c(\bar \nabla)+d\tilde \hA^c(\tilde \nabla^{T^v\pi}\oplus\tilde  \nabla^{\underline{\R^N}_E },\bar \nabla_{|\{t\}\times E})\\
&=&-dt\wedge i_{\partial_t} \hA^c(\bar \nabla)
 -\hA^c(\bar \nabla_{|\{t\}\times E})+\hA^c(\tilde \nabla^{T^v\pi})\ .
\end{eqnarray*}
We get
$$ R(\bar o)=\pr_E^*\hA^c(\tilde \nabla^{T^v\pi})-d\bar \sigma=
\hA^c(\bar \nabla)\ .$$
In combination with (\ref{ui1}) and (\ref{ui2}) this implies 
$\hat \rho^k(o_\pi)= 1$.
\hB

\section{The index theorem}

Let $\pi:E\to B$ be a proper submersion over a compact base $B$ with fibre dimension $n:=\dim(E)-\dim(B)$. We assume that $\pi$ is topologically $K$-oriented by the datum of a $Spin^c$-structure on $T^v\pi$. Let $o$ be a smooth $K$-orientation of $\pi$ which refines this topological $K$-orientation. Then we have the push-forward
$$\hat \pi_!:\hat K^*(E)\to \hat K^{*-n}(B)\ ,$$
see  \cite[Definition 3.18]{bunke-20071}.  The following theorem refines the identity (\ref{eq20})
to the smooth case.

 \begin{theorem}\label{uifqfefewf5454577}
In $\hat K^{*-n}(B)[\frac{1}{k}]$  we have the identity
\begin{equation}\label{uqdqwdwqdwqdqd31}
\hat \Psi^k(\hat \pi_!(\hat x))= \hat \pi_!(\hat \rho^k(o)^{-1}\cup \hat \Psi^k(\hat x))\ , \quad \forall x\in \hat K^*(E) [\frac{1}{k}]
\end{equation}
\end{theorem}
\proof  

\begin{lem}\label{uieewfwefwef}
The equality (\ref{uqdqwdwqdwqdqd31}) holds true after applying $I$ or $R$. Moreover
it holds true if $\hat x=a(\alpha)$ for $\alpha\in \Omega P^{*-1}(E)$.
\end{lem}
\proof
The proof goes by straightforward calculations.
\hB
 
We now consider the difference of the left- and right-hand sides of (\ref{uqdqwdwqdwqdqd31}):
$$\hat \Delta_\pi(\hat x):=\hat \Psi^k(\hat \pi_!(\hat x))-\hat \pi_!(\hat \rho^k(o)\cup \hat \Psi^k(\hat x))\in \hat K^{*-n}(B)[\frac{1}{k}]\ .$$
By Lemma \ref{uieewfwefwef} we know that $R\circ\hat \Delta_\pi=0$ and $I\circ \hat \Delta_\pi=0$. It follows that
$$\hat \Delta_\pi(\hat x)\in
HP^{*-n-1}(B)/\im(\ch)[\frac{1}{k}]\ .$$
Moreover, since  it vanishes on classes of the form $\hat x=a(\alpha)$, it factors over a 
 homomorphism
$$\Delta_\pi:K^*(E)\to HP^{*-n-1}(B)/\im(\ch)[\frac{1}{k}]\ .$$

\begin{lem}\label{zwqudqwdqwd}
Assume that $p:F\to B$ is a smoothly $K$-oriented zero bordism of the smoothly $K$-oriented  bundle $\pi:E\to B$, and that $y\in K^l(F)[\frac{1}{k}]$.  Then we have $\Delta_\pi(y_{|E})=0$.
\end{lem}
\proof
We let $o_p$ denote a smooth $K$-orientation of $p$ with a product structure near the boundary which restricts to a  smooth $K$-orientation $o$ of $\pi$.  
We further choose a smooth lift $\hat y\in \hat K^l(F)$. We let $\hat x:=\hat y_{|E}$ and $x:=I(\hat x)=y_{|E}$.
Then we calculate using the bordism formula \cite[Prop. 5.18]{bunke-20071}
\begin{eqnarray*}
\Delta_\pi(x)&=&\hat \Psi^k(\hat \pi_!(\hat x))- \hat \pi_!(\hat \rho^k(o)\cup \hat\Psi^k(\hat x))\\&=&\hat \Psi^k\left( a\left(\int_{F/B}R(o_p)\wedge R(\hat y)\right)\right)-a \left(\int_{F/B} R(o_p) \wedge R(\hat\rho^k(o_p))\wedge R(\hat \Psi^k(\hat y))\right)\\
&=& a \left(\Psi^k_\Omega\left(  \int_{F/B} R(o_p)\wedge R(\hat y)\right) -\int_{F/B}  \Psi^k_\Omega (R(o_p)) \wedge  \Psi^k_\Omega (R(\hat y))\right) \\
&=&  0
\end{eqnarray*}
\hB

\begin{lem}\label{zwqudqwdqwd1}
The homomorphism $\Delta_\pi$ only depends on the underlying 
topological $K$-orientation of the  bundle $\pi:E\to B$.
\end{lem}
\proof
Let $o_0$ and $o_1$ be two smooth $K$-orientations with the same underlying
topological $K$-orientations which gives  rise to $\Delta_{\pi,0}$ and $\Delta_{\pi,1}$.
Then we choose a smooth $K$-orientation $o_p$ on $p:=\pi\circ\pr_E:F:=[0,1]\times E \to B$ which restricts to  $o_i$ on the endpoints of the interval.
We apply Lemma \ref{zwqudqwdqwd}
to the class
$y=\pr_E^*x$ in order to see that
$\Delta_{\pi,0}(x)=\Delta_{\pi,1}(x)$. 
\hB

The homomorphism $\Delta_\pi$ has the following naturality property. Let
 $$\xymatrix{E^\prime\ar[d]^{\pi^\prime}\ar[r]^G&E\ar[d]^\pi\\B^\prime\ar[r]^g&B}$$
be  cartesian with a compact manifold $B^\prime$ and the topological $K$-orientation of $\pi^\prime$ be  induced by that of $\pi$.
\begin{lem}\label{uieowfwefewfwefw} We have 
$$\Delta_{\pi^\prime}\circ G^*=g^*\circ \Delta_\pi\ .$$
\end{lem}
\proof
This follows from
the naturality of $\hat \Psi^k$, $G^*\hat \rho^k(o_\pi)=\hat \rho^k(o_{\pi^\prime})$ (Lemma \ref{uiqwdwqdwqdwqdwqd})
and
$\hat \pi_!\circ G=g^*\circ \hat \pi_!$ (see \cite[Lemma 3.20]{bunke-20071}).
\hB 
 
The following proposition is the nontrivial heart of the proof of Theorem \ref{uifqfefewf5454577}.
\begin{prop}\label{dzqwudw1}
If $\dim(E)=n=2m-1$ and $B=*$, then $\Delta_\pi=0$.
\end{prop}
\proof
Let $E$ be a closed $Spin^c$-manifold of  dimension $2m-1$ together with a class $x\in K^0(E)$  which we view as a homotopy class of maps
$E\to \Z\times BU$. The pair $(E,x)$ represents a $Spin^c$-bordism class
$[E,x]\in \Omega^{Spin^c}_{2m-1}(\Z\times BU)$.
The integral cohomology of $\Z\times BU$ is concentrated in even degrees, 
and the  odd part of $\Omega^{Spin^c}_*$ is a torsion $2$-group.
 Using the  Atiyah-Hirzebruch spectral sequence we see that  $\Omega^{Spin^c}_{2m-1}(\Z\times BU)$
is a torsion $2$-group. Hence there exists $l\in \nat$ of the form $l=2^{l^\prime}$ such that
$l[E,x]=0$. Thus there exists a $2m$-dimensional $Spin^c$-manifold
$W$  with boundary
$\partial W\cong l E$ together with an extension $y:W\to \Z\times BU$ of the map $
\partial W\to \Z\times BU$ induced by $x$.
More precisely, the boundary of $\partial W$ decomposes as
$\partial W=\bigsqcup_{i=1}^l\partial_i W$, and we can choose identifications of $Spin^c$-manifolds $w_i:E\stackrel{\sim}{\to} \partial_i W$ such that
$x=y\circ w_i$ for all $i=1,\dots,l$. 
We choose a $\Z/2\Z$-graded vector bundle $Y\to W$ such that its $K$-theory class 
satisfies $[Y]=y$. We further define $X:=w_1^*Y$ such that $[X]=x$. After stabilisation of $Y$, if necessary, we can assume that there are isomorphisms
$w_i^*Y\cong X$ for all $i=1,\dots,l$.
These choices constitute a $Spin^c$-$\Z/l\Z$-manifold $\check W$ in the sense of \cite{MR1144425} and  \cite{MR1060635} together with a $\Z/l\Z$-bundle $\check Y$ over $\check W$.

We consider $\R^2\cong \C$. Let $\xi$ be a primitive  $l$'th root of unity.
We fix $r>0$ so small that the discs $B(\xi^i,r)$ are pairwise disjoint. We let $\check \R^2 :=\R^2\setminus \bigcup_{i=0}^{l-1}\inter B(\xi^i,r)$. In order to define the structure of a $\Z/l\Z$-manifold $\check \R^2$
we fix the identifications
$v_i:S^1\to  \partial_i\check \R^2=\partial B(\xi^i,r)$ as
$v_i(u)=\xi^i+u r\xi^i$.
Let
$\bar \R^2$
be the quotient of $\check \R^2 $ obtained by identifying the boundary components with $S^1$ using the maps $v_i$. 

As shown in \cite{MR1144425} we have $\tilde K^0(\bar \R^2)\cong \Z/l\Z$.
Let us describe this isomorphism explicitly.
We choose $R>4$ and consider the  decomposition
$\bar \R^2\cong (\R^2\setminus \inter B(0,R))\cup_{S(0,R)} B(0,R)\cap \bar \R^2$.
We define line bundles $L_m$ on $\bar \R^2$ using a clutching function
$S(0,R)\cong S^1\to U(1)$ of degree $m$.  In greater detail, we define
$L_m:= ((\R^2\setminus \inter B(0,R))\times \C \sqcup B(0,R)\cap \bar \R^2\times \C)/\sim$,
where the glueing is given by 
$(u^\prime,z)\sim (u, (\frac{u}{R})^mz)$
if $u^\prime\in S(0,R)\cap (\R^2\setminus \inter B(0,R))$,
$u\in S(0,R)\cap (B(0,R)\cap \bar \R^2)$,
$u=u^\prime$ as points in $\R^2=\C$.

If $m$ is divisible by 
$l$, then the clutching function   extends to $B(0,R)\cap \bar \R^2$ and therefore defines the trivial line bundle. The element $[m]\in \Z/l\Z\cong \tilde K^0(\bar \R^2)$ is represented by
$L_m-L_0$.
With this description it is easy to see how the Adams operation $\Psi^k$ acts on $K^0(\bar \R^2)$.
Indeed, $\Psi^k([L_m-L_0])=[L_m^k-L_0^k]=[L_{km}-L_0]$.
Note that by Bott periodicity $\tilde K^{2-2g}(\bar \R^2)\cong \Z/l\Z$, too.
Therefore  on $\tilde K^{2-2g}(\bar \R^2)\cong \Z/l\Z$ the action of the Adams operation $\Psi^k$ is given by multiplication by 
$k^{g}$.

We now recall the definition of  the topological $\Z/l\Z$-index  given by \cite{MR1060635}. 

We choose collars $c_i:[0,1)\times E\to W$.
Then we define the map of $\Z/l\Z$-manifolds
$\check \rho:\check W\to \check  \R^2 $ as follows.
On the $i$th collar we require that
$$\xymatrix{[0,1]\times E\ar[d]^{\pr_1}\ar[r]^{c_i}& \check W\ar[d]^{\check \rho}\\[0,1]\ar[r]^{\gamma_i}&\check\R^2}$$ commutes, where $\gamma_i: [0,1]\to \R^2$ is given by 
$\gamma_i(t):=(1-t)(1-r)\xi^i$.
The complement of the union of collars is mapped to the origin $0\in  \check \R^2 $.

Let furthermore $i:W\to V$ be an embedding of the manifold $W$ into a real vector space $V$. Then we consider the embedding of
$\Z/l\Z$-manifolds
$i\times \check \rho:\check W\to V\times \check \R^2$.
We can choose a $\Z/l\Z$-normal bundle
$\check N\to \check W$ and extend the embedding to an open embedding
$\check I:\check N\to V\times \check \R^2$.
We let
$\bar N\to \bar W$ be the quotient obtained by identifying the boundary components to one.
We get an induced map
$\bar I:\bar N\to V\times \bar \R^2$.
The bundle $\bar N\to \bar W$ has an induced $Spin^c$-structure and therefore has a  Thom isomorphism $\Thom_{\bar N}$.
The topological index $\ind^{\Z/l\Z}_t:K^0(\bar W)\to \Z/l\Z$ is given by
\begin{equation}\label{uiwqdqwdwqdqwdq3232}
\ind^{\Z/l\Z}_t:K^*(\bar W)\stackrel{\Thom_{\bar N}}{\to}K_c^{*+\dim(\bar N)}(\bar N)\stackrel{excision}{\to} K^{*+\dim(\bar N)}_c(V\times \bar \R^2)\stackrel{Thom^{-1}}{\to}
\tilde K^{*-\dim(W)}(\bar \R^2)\ .
\end{equation}

We now choose a smooth $K$-orientation $o_E$ of the $K$-oriented map $\pi:E\to *$ and a geometry $\bX$ on $X$. We extend this orientation to an orientation $o_W$ of the $W\to *$. Similarly we extend the geometry $\bX$ to a geometry $\bY$ of $Y$. In this way we get geometric manifolds
$\cW$ and $\cE$ such that $\partial \cW= l\cE$ (see \cite[Def. 2.1.30]{bunke-2002}). 
We let $\cE\otimes \bX$ denote the geometric manifold  obtained from $\cE$ 
by twisting the Dirac bundle of $\cE$ with $\bX$.
Since $\ind(\cE\otimes \bX)=0$ (e.g. since $l(\cE\otimes \bX)$ is zero bordant)
we can  choose a taming $(\cE\otimes \bX)_t$. It gives a boundary taming $(\cW\otimes \bY)_{bt}$.
The class
$$\ind_a^{\Z/l\Z}(D(\cW\otimes \bY)):=[\ind(\cW\otimes \bY)_{bt}]\in \Z/l\Z$$
is the analytic $\Z/l\Z$-index. We refer to \cite[Def. 2.1.44, Def. 2.1.47]{bunke-2002} for the definition of a taming or  boundary taming and the corresponding index theory.

The index theorem of \cite{MR1144425},  \cite{MR1060635} states that
\begin{equation}\label{uidwqdwqdwdd65656556}
\ind_a^{\Z/l\Z}(D(\cW\otimes \bY))=\ind_t^{\Z/l\Z}(\bar y)\ ,
\end{equation}
where $\bar y\in K^0(\bar W)$ is represented by the map
$\bar y:\bar W\to \Z\times BU$ induced by $y$, and on the right-hand side we use the identification $\tilde K^{*-2m}(\bar \R^2)\cong \Z/l\Z$ given above.
We apply the $APS$-index theorem for boundary tamed manifolds \cite[Thm.2.2.18]{bunke-2002}
and get
\begin{equation}\label{uiedqwdwqdwqq89899898}
\ind(\cW\otimes \bY)_{bt}=\Omega(\cW\otimes \bY)-l\eta((\cE\otimes \bX)_{t})\ ,
\end{equation}
where the first term is the usual local contribution and the second term is the boundary correction.
We now choose a finite  formal sum $\bQ:=\sum_\alpha a_\alpha \bQ_\alpha$ of geometric bundles $\bQ_\alpha\to W$ with coefficients in $a_\alpha \in \Z[\frac{1}{k}]$ 
 which represents $\hat \rho^k(o_E)\in \hat K^0(W)[\frac{1}{k}]$.  This is possible, see e.g. \cite{simons-2008}).
More precisely, if $\cQ_\alpha$ denotes the geometric family induced by $\bQ_\alpha$, then
we assume that $\sum_\alpha a_\alpha[\cQ_\alpha,0]=\hat \rho^k(o_E)$.
In the following we will suppress this sum decomposition.
The pull-back $\bR:=c_1^*\bQ$ represents $\hat \rho^k(o_E)\in \hat K^0(E)$.
Indeed we can choose the geometry of $\bQ$ such that we have isomorphisms
$c_i^*\bQ\cong c_1^*\bQ$. In this way we get a $\Z/l\Z$-bundle $\check \bQ$ over $\check W$.
Its underlying topological $\Z/l\Z$-bundle represents $\rho^k(\bar N)$.
From the construction of the topological index (\ref{uiwqdqwdwqdqwdq3232})  and  the calculation of the action of the Adams operation on $\tilde K^{2-2m}(\bar \R^2)$ given above (note that
$\ind_t^{\Z/l\Z}(\bar y)\in \tilde K^{2-2m}(\bar \R^2)$)
we have the following indentities
\begin{eqnarray*}
\Psi^k\ind_t^{\Z/l\Z}(\bar y)&=&\ind_t^{\Z/l\Z}(\rho^k(\bar N)\cup \bar y)\\
\Psi^k\ind_t^{\Z/l\Z}(\bar y)&=&k^{m} \ind_t^{\Z/l\Z}(\bar y)\ .
\end{eqnarray*}
We let $[\bX]\in \hat K^0(E)$ denote the smooth $K$-theory class induced by the geometric bundle $\bX$. The following calculation uses the explicit cycle level description of the push-forward in smooth $K$-theory \cite[(17)]{bunke-20071}  and the relations  \cite[Def. 2.10]{bunke-20071}. We get
\begin{eqnarray*}
\Delta_\pi(x)&=&\hat \Psi^k(\hat \pi_!([\bX]))-\hat \pi_!( \hat \rho^k(o_E)\cup[\bX])\\
&=&\hat \Psi^k([\cE\otimes \bX,0])-[\bE\otimes \bX\otimes \bR,0])\\
&=&\Psi^k([\emptyset,\eta((\cE\otimes \bX)_{t})])-[\emptyset,\eta((\bE\otimes \bX\otimes \bR)_{t})]\\
&=&a\left(k^{m}\eta((\cE\otimes \bX)_{t})-\eta((\bE\otimes \bX\otimes \bR)_{t})\right)
\end{eqnarray*}
In $\R/\Z$ we have by (\ref{uiedqwdwqdwqq89899898}) the following identity 
\begin{eqnarray*}
\lefteqn{[k^{m}\eta((\cE\otimes \bX)_{t})-\eta((\bE\otimes \bX\otimes \bR)_{t})]_{\R/\Z}}&&\\
&=&[k^{m} l^{-1}\Omega(\cW\otimes \bY)-
l^{-1}\Omega(\cW\otimes \bY\otimes \bQ)]_{\R/\Z}\\&&+l^{-1}\ind_a^{\Z/l\Z}(D(\cW\otimes \bY\otimes \bQ))-l^{-1}k^{m}\ind_a^{\Z/l\Z}(D(\cW\otimes \bY))\ ,
\end{eqnarray*}
where we interpret $\Z/l\Z\subset \R/\Z$ via multiplication by $l^{-1}$.
 We now observe that
$$k^{m} \Omega(\cW\otimes \bY)=\Omega(\cW\otimes \bY\otimes \bQ)$$
and that  in $\R/\Z[\frac{1}{k}]$  we have by (\ref{uidwqdwqdwdd65656556}) the identity
\begin{eqnarray*}\lefteqn{
k^{m}[l^{-1}\ind_a^{\Z/l\Z}(D(\cW\otimes \bY)) -
[l^{-1}\ind_a^{\Z/l\Z}(D(\cW\otimes \bY\otimes \bQ))] }&&\\&=&
l^{-1}k^{m} \ind_t^{\Z/l\Z}(\bar y)-l^{-1}\ind_t^{\Z/l\Z}( \rho^k(\bar N)\cup\bar y)\\
&=&0\ .
 \end{eqnarray*}
This implies that in $\R/\Z[\frac{1}{k}]$ we have 
$$k^{m}\eta((\cE\otimes \bX)_{t})-\eta((\bE\otimes \bX\otimes \bR)_{t})=0$$
and hence $\Delta_\pi(x)=0$. \hB

\begin{lem}\label{uiwqdwqdqwdwdqwd} Let  $n:=\dim(E)-\dim(B)$ be even  and $x\in K^0(E)$. Then
 $$\Delta_\pi(x)=0\in HP^{-n-1}(B)/\im(\ch)[\frac{1}{k}]\ .$$
\end{lem}
\proof
For $x\in K^{0}(E)$ we have
$$\Delta_\pi(x)\in HP^{-n-1}(B)/\im(\ch)[\frac{1}{k}]\subseteq \hat K_{flat}^{-n}(B)[\frac{1}{k}]\cong K\R/\Z^{-n-1}(B)[\frac{1}{k}]$$
(see \cite[Thm 5.5]{bs2009} for the last isomorphism).
We use the universal coefficient formula \cite{MR0388375}
$$K\R/\Z^{-n-1}(B)[\frac{1}{k}]\cong \Hom(K_{-n-1}(B),\R/\Z)[\frac{1}{k}]\ .$$
In order to show that
$\Delta_\pi(x)=0$ it therefore suffices to show that
$$\langle u,\Delta_\pi(x)\rangle =0\in\R/\Z[\frac{1}{k}] $$
for all $K$-homology classes $u\in K_{-n-1}(B)$.
We now use the geometric picture of $K$-homology \cite{MR679698}, \cite{MR2330153}.
Given $u\in K_{-n-1}(B)$ there exists a $k$-dimensional $Spin^c$-manifold $Z$ (where $k$ is odd) together with a map
$f:Z\to B$ such that
$u=b^jf_*( [Z])$. Here $[Z]\in K_{k}(Z)$ is the $K$-homology orientation of $Z$ given by the $Spin^c$-structure and $j:=\frac{-n-1-k}{2}$.
Let $q:M\to *$ be the projection.
Then we have for $z\in K\R/\Z^{-n-1}(B)$
$$\langle u,z\rangle=q_!(  b^{-j}f^*z)\in K\R/\Z^0\cong \R/\Z\ .$$
We consider the diagram
$$\xymatrix{W\ar@/_1cm/@{.>}[dd]_r \ar[d]^p\ar[r]^g&E\ar[d]^\pi \\ Z\ar[d]^q\ar[r]^f&B\\{*}&}\ .$$
We choose smooth $K$-orientations $o_q$ and $o_\pi$ on $q$ and $\pi$ lifting the topological ones.
Furthermore we choose a smooth lift $\hat x$   of $x$. We equip $p$ with the induced smooth $K$-orientation $o_p$. By Lemma \ref{suskwswqoidd} there exists a class $\hat z\in \hat K^0(Z)[\frac{1}{k}]$ such that $\hat \Psi^k(\hat z)=\hat \rho^k(o_q)^{-1}$.
Then we calculate using Proposition \ref{dzqwudw1} at the marked equalities,  omitting the powers of the Bott element $b^j$ in order to increase readability, and using the projection formula  \cite[Prop. 4.5]{bunke-20071}
\begin{eqnarray*}
\lefteqn{\langle u,\Delta_\pi(x)\rangle}&&\\&=&
q_!( f^*\Delta_\pi(x))\\&=&
\hat q_!( \ f^*\hat \Delta_\pi (\hat x))\\
&\stackrel{Lemma \ref{uieowfwefewfwefw}}{=}&q_!(  \hat \Delta_p (g^*\hat x))\\
&=&\hat q_!( \hat \Psi^k(\hat p_!(\hat x))-\hat p_!(\hat \rho^k(o_p)\cup \hat \Psi^k(\hat x)))\\
&=&\hat q_!(\hat \rho^k(o_q) \cup \hat \rho^k(o_q)^{-1}\cup \hat \Psi^k(\hat p_!(\hat x))-\hat p_!(p^*\hat \rho^k(o_q)\cup p^*\hat \rho^k(o_q)^{-1}\cup \hat \rho^k(o_p)\cup \hat \Psi^k(\hat x)))\\
&=&\hat q_!(\hat \rho^k(o_q) \cup  \hat \Psi^k(\hat z)\cup \hat \Psi^k(\hat p_!(\hat x))-\hat p_!(p^*\hat \rho^k(o_q)\cup  p^* \hat \Psi^k(\hat z)\cup \hat \rho^k(o_p)\cup \hat \Psi^k(\hat x)))\\
&=&\hat q_!(\hat \rho^k(o_q) \cup  \hat \Psi^k(\hat z \cup  \hat p_!(\hat x)))-\hat q_!\circ \hat p_!(p^*\hat \rho^k(o_q)\cup     \hat \rho^k(o_p)\cup \hat \Psi^k(p^*\hat z\cup \hat x))\\
&\stackrel{!}{=}&\hat \Psi^k (\hat q_!(\hat z\cup \hat p_!(\hat x)))-\hat q_!\circ \hat p_!( \hat \rho^k(o_r)\cup \hat \Psi^k(p^*\hat z\cup \hat x))\\
&=&\hat \Psi^k (\hat r_!( p^*\hat z\cup\hat x))-\hat r_!(\hat \rho^k(o_r)\cup  \hat \Psi^k(p^* \hat z\cup \hat x))\\
&\stackrel{!}{=}&0\ .
\end{eqnarray*}
\hB 
\begin{lem}\label{uiwqdwqdqwdwdqwd1}
Let $n:=\dim(E)-\dim(B)$ be odd   and $x\in K^{-1}(E)$. Then
 $$\Delta_\pi(x)=0\in H P^{-n-1}(B)/\im(\ch)[\frac{1}{k}]\ .$$
\end{lem}
\proof
We consider the bundle $q=\pi\circ \pr_E:S^1\times E\to B$ with even-dimensional fibres
and the class $y=e\times x\in K^0(S^1\times E)$, where $e\in K^1(S^1)\cong \Z$ is the generator.
We choose smooth lifts $\hat e$ and $\hat x$ of $e$ and $x$. Furthermore we choose a smooth $K$-orientation $o_p$ lifting the underlying topological $K$-orientation, and we  let $o_{\pr_E}$ be such that
$\int=\pr_{E,!}$.
Then we have $\pr_{E,!}(\hat e\times \hat x)=\hat x$, $\hat \Psi^k(\hat e)=\hat e$, and
$\pr_E^*\hat \rho^k(o_\pi)=\hat \rho^k(o_q)$. Applying Lemma \ref{uiwqdwqdqwdwdqwd} to
$\Delta_q(e\times x)$ we get
\begin{eqnarray*}
\Delta_\pi(x)&=&\Psi^k(\hat \pi_!(\hat x))-\hat \pi_!(\hat \rho^k(o_\pi)\cup \hat \Psi^k(\hat x))\\
&=&\hat \Psi^k(\hat q_!(\hat e\times  \hat x))-\hat q_!(\pr_E^*\hat \rho^k(o_{\pi})\cup \hat e\cup \pr_E^*\hat \Psi^k(\hat x))\\&=&\Delta_q(e\times x)\\&=&0\ .
\end{eqnarray*}
\hB

\begin{lem}\label{uiwqdwqdqwdwdqwd2} Let $n:=\dim(E)-\dim(B)$ be even   and
  $x\in K^{-1}(E)$. Then
 $$\Delta_\pi(x)=0\in HP^{-n-1}(B)/\im(\ch)[\frac{1}{k}]\ .$$
\end{lem}
\proof
We consider the diagram
$$\xymatrix{S^1\times E\ar[r]^{\pr_E}\ar[d]_{q:=\id\times \pi}&E\ar[d]^\pi\\
S^1\times B\ar[r]^{\pr_B}&B}\ .$$
We claim that
$$e\times \pr_B^*\Delta_\pi(x)=e\times \Delta_q(\pr_E^*x)\ .$$
Indeed, after choosing
 smooth orientations and smooth lifts we calculate using
$\hat \rho^k(o_q)=\pr_E^*\hat \rho^k(o_p)$ (Lemma \ref{uiqwdwqdwqdwqdwqd}), $\hat \Psi^k(\hat e)=\hat e$, the equality
$\hat e\times \hat \Psi(\hat y)=\hat \Psi^k(\hat e\times \hat y)$,
and $\hat q_!(\hat e\times \hat y)=\hat e\times \hat \pi_!(\hat y)$ (a special case of the projection formula \cite[Prop. 4.5]{bunke-20071}) 
\begin{eqnarray*}
e\times \Delta_q(\pr_E^*x)&=&\hat e\times (\hat \Psi^k(\hat \pi_!(\hat x))-\hat \pi_!(\hat \rho^k(o_\pi)\cup \hat \Psi^k(\hat x)))\\
&=& \hat \Psi^k(\hat e\times \hat \pi_!(\hat x))-\hat q_!(\hat e\times (\hat \rho^k(o_\pi)\cup \hat \Psi^k(\hat x)))\\
&=&\hat \Psi^k(\hat q_!(\hat e\times \hat x))-\hat q_!(\hat \rho^k(o_q)\cup \hat \Psi^k(\hat e\times \hat x))\\
&=&\Delta_q(e\times x)\\
&\stackrel{Lemma \:\ref{uiwqdwqdqwdwdqwd}}{=}&0\ .
\end{eqnarray*} \hB

\begin{lem}\label{uiwqdwqdqwdwdqwd3} Let  $n:=\dim(E)-\dim(B)$ be odd   and   $x\in K^{0}(E)$. Then
 $$\Delta_\pi(x)=0\in HP^{-n-1}(B)/\im(\ch)[\frac{1}{k}]\ .$$
\end{lem}
\proof
Let $q:=\pi\circ \pr_E:S^1\times E\to B$. This bundle has even-dimensional fibres. 
We calculate (again after choosing smooth lifts $\hat x$ and $\hat e$ of $x$ and $e$ and smooth $K$-orientations  $o_\pi$, $o_{\pr_E}$ refining the underlying topological ones)
\begin{eqnarray*}
\Delta_\pi(x)&=&\hat \Psi^k(\hat \pi_!(\hat x))-\hat \pi_!(\hat \rho^k(o_\pi)\cup \hat \Psi^k(\hat x))\\
&=&\hat \Psi^k(\hat q_!(\hat e\times \hat x))-\hat q_!(\hat e\times \hat \rho^k(o_\pi)\cup \hat \Psi^k(\hat x))\\
&=&\hat \Psi^k(\hat q_!(\hat e\times \hat x))-\hat q_!(  \hat \rho^k(o_q)\cup \hat \Psi^k(\hat e\times \hat x))\\
&=&\Delta_q(e\times x)=0\ .
\end{eqnarray*}
 \hB 

The collection of the Lemmas \ref{uiwqdwqdqwdwdqwd}, \ref{uiwqdwqdqwdwdqwd1},
\ref{uiwqdwqdqwdwdqwd2} and \ref{uiwqdwqdqwdwdqwd3} gives the Theorem \ref{uifqfefewf5454577}
\hB

\section{Application to the $e$-invariant}\label{sec4}

Let $\pi:E\to B$ be a proper submersion over a compact base with fibre-dimension $n$ together with a stable framing of $T^v\pi$. In this situation we have the canonical smooth $K$-orientation $o_\pi$ and the class
$e(\pi):=\hat \pi_!(1)\in \hat K^{-n}(B)$, see  \cite[5.11]{bunke-20071}. 
This class is actually flat and therefore belongs to
$\hat K^{-n}_{flat}(B)\cong K\R/\Z^{-n-1}(B)$. It is an invariant of the bordism class
of bundles with stably framed vertical bundles.

\begin{theorem}\label{udqiwdqwdqwd}
The $e$-invariant satisfies
$$(\hat \Psi^k-1)e(\pi)=0\in  K\R/\Z^{-n-1}(B)[\frac{1}{k}]\ .$$
\end{theorem}
\proof
Since $\hat \rho^k(o_\pi)=1$ by Proposition \ref{udiqwdwqdw33333}
using Theorem \ref{uifqfefewf5454577} we get
$$\hat \Psi^k e(\pi)=\hat \Psi^k(\hat \pi_!(1))=\hat \pi_!(\hat \rho^k(o_\pi))=\hat \pi_!(1)=e(\pi)\ .$$
\hB

In the special case that
$B=*$ and $n=2m-1$ is odd, 
$e(\pi)\in K\R/\Z^{-n-1}\cong \R/\Z$ is the $e$-invariant of Adams of the framed bordism class
$[E]\in \Omega^{fr}_n$ represented by $E$, or equivalently of the element in the stable stem $\pi^S_{n}$ corresponding to $[E]$ via the Pontrjagin-Thom construction. 
Originally, the $e$-invariant has been introduced in order to detect elements in the image of the $j$-homomorphism. The order of the image of $j$-homomorphism is known (see the series of papers
 \cite{MR0198468} and conforms with the following special case of Theorem \ref{udqiwdqwdqwd}.
\begin{kor}\label{udidqwdwd}
If $\pi:E\to *$ is the projection from a compact stably framed manifold to the point and $\dim(E)=2m-1$, then
 for every $k\in \{-1\}\cup \nat$ there exists $L\in \nat$ such that
$$k^L(k^{m}-1)e(\pi)=0\ .$$
\end{kor}
\proof
Indeed
$0=(\hat \Psi^k-1)e(\pi)=(k^m-1)e(\pi)\in \R/\Z[\frac{1}{k}]$.
\hB 
The determination of the order of the image of the $j$-homomorphism in the work of Adams
also uses Adams operations, namely in order to characterise the kernel of $j$. The proof
of the upper bound of the order of the $e$-invariant Corollary \ref{udidqwdwd} seems to employ the Adams operations in a completely different manner.

 Let us now discuss some application to higher $\rho$-invariants. Let $\pi:E\to B$ be a proper submersion over a compact base $B$ which is $K$-oriented by a $Spin^c$-structure on the vertical bundle $T^v\pi$. We fix a base point of $E$, chose a character $\chi:\pi_1(E,*)\to \Z/k\Z$, and we let $\bH:=(H,h^H,\nabla^H)$ be the corresponding
geometric line bundle. It represents a class $[\bH]\in \hat K^0(E)$. Since
$\bH^k$ is trivial we have the relation $[\bH]^k=1$ in $\hat K^0(E)$.
We choose  a smooth $K$-orientation $o_\pi$ which refines the topological one. The higher $\rho$-invariant is then defined by  
$$\rho(\chi):=\hat \pi_!([\bH])-\hat \pi_!(1)\ .$$
Note that  $R(\rho(\chi))=0$ so that $\rho(\chi)\in K^{-n}_{flat}(B)\cong K\R/\Z^{-n-1}(B)$, and it is independent of the choice of $o_\pi$. Therefore $\rho(\chi)$ is a differential-topological invariant of the $K$-oriented bundle $\pi$.

\begin{prop}
We have the relation
$\hat \Psi^k(\rho(\chi))=0\in  K\R/\Z^{-n-1}(B)[\frac{1}{k}]$.
\end{prop}
\proof
We calculate  
$$\hat \Psi^k(\hat \pi_!(1))-\hat \Psi^k(\hat \pi_!([\bH]))=\hat \pi_!(\hat \rho^k(o_\pi))-\hat \pi_!(\hat \rho^k(o_\pi)\cup \hat \Psi^k([\bH]))= 0\in K\R/\Z^{-n-1}(B)[\frac{1}{k}]$$ 
since
$\hat \Psi^k([\bH])=[\bH]^k=1$
by Proposition \ref{udiduqwdwqddwqdwqdw} .\hB
If $M$ is a point and $\dim(E)=2m-1$, then
$\rho(\chi)\in K\R/\Z^{-2m}(*)\cong \R/\Z$. 
Furthermore
$\hat \Psi^k(\rho(\chi))=k^{m} \rho(\chi)$.
\begin{kor}
There exists $L\in \nat$ such that
$$k^L\rho(\chi)=0\in \R/\Z\ .$$
\end{kor}
This can also be verified independently from the present formalism.



\end{document}